\documentclass[11pt,a4paper]{article}
\usepackage[a4paper, total={6in, 8in}]{geometry}
\usepackage{authblk}

\usepackage[utf8]{inputenc}
\usepackage{hyperref}
\usepackage{tabularx}
\usepackage{enumitem} 
\newlist{condenum}{enumerate}{1} 
\setlist[condenum]{label=\bfseries P\arabic*., 
                   ref=\arabic*, wide}

\usepackage{algorithm}

\usepackage{algpseudocode}
\usepackage[cmex10]{amsmath}
\usepackage{%
	amsfonts,%
	amsmath,%
	amssymb,%
	amsthm,%
	bbm,%
	cite,%
	float,%
	etoolbox,%
	mathrsfs,%
	mathtools,%
	multirow,%
	pgf,%
	pgfplots,%
	subfigure,%
	tikz,%
}
\usepackage[normalem]{ulem}
\usepackage{cleveref}
\pgfdeclareplotmark{mystar}{
    \node[star,star point ratio=2.25,minimum size=6pt,
          inner sep=0pt,draw=black,solid,fill=red] {};
}
\pgfdeclareplotmark{mystar2}{
    \node[star,star point ratio=2.25,minimum size=6pt,
          inner sep=0pt,draw=black,solid,fill=blue] {};
}

\newlist{temp}{enumerate}{1} 
\setlist[temp]{label=(\roman*), 
                   ref=\roman*, wide}
\usepgflibrary{shapes}
\usetikzlibrary{%
	arrows,%
	patterns,%
	positioning,%
	shapes.callouts,%
	shapes.symbols,%
}
\usepackage{bbm}
\usepackage{breqn}

\newcommand\numeq[1]%
  {\stackrel{\scriptscriptstyle(\mkern-1.5mu#1\mkern-1.5mu)}{=}}

\theoremstyle{plain}
\newtheorem{thm}{Theorem}
\newtheorem{lem}[thm]{Lemma}
\newtheorem{proposition}[thm]{Proposition}

\theoremstyle{definition}

\newtheoremstyle{case}{}{}{}{}{}{:}{ }{}
\theoremstyle{case}

\theoremstyle{remark}

\newcommand{\eq}[1]{\begin{align*}#1\end{align*}}
\newcommand{\EQ}[1]{\begin{equation*}#1\end{equation*}}
\newcommand{\eqn}[1]{\begin{align}#1\end{align}}
\newcommand{\EQN}[1]{\begin{equation}#1\end{equation}}

\newcommand{\norm}[1]{\left\lVert#1\right\rVert}

\newcommand{\E}{\mathbb{E}}
\newcommand{\Z}{\mathbb{Z}}

\newcommand{\mb}[1]{\mathbb{#1}}
\newcommand{\mc}[1]{\mathcal{#1}}
\newcommand{\mf}[1]{\mathbf{#1}}

\newcommand{\brac}[1]{\left(#1\right)}
\newcommand{\cbrac}[1]{\left\{#1\right\}}
\newcommand{\sbrac}[1]{\left[#1\right]}
\newcommand{\indic}[1]{\mathbbm{1}{\brac{#1}}}

\newcommand{\define}{\triangleq}

\title{The Power of Two Choices with Load Comparison Errors}
\author[$*$]{Sanidhay Bhambay}
\author[$*$]{Arpan Mukhopadhyay}
\author[$\dag$]{Thirupathaiah Vasantam}
\affil[$*$]{Department of Computer Science, University of Warwick}
\affil[$\dag$]{Department of Computer Science, Durham University}

\date{}

\begin{document}
\maketitle
\begin{abstract}
We consider a system with $n$ unit-rate servers where jobs arrive according a Poisson process with rate $n\lambda$ ($\lambda <1$). In the standard \textit{Power-of-two} or Po2 scheme, for each incoming job, a job dispatcher samples two servers uniformly at random and sends the incoming job to the least loaded of the two sampled servers. 
However, in practice, the load information may not be accurate at the job dispatcher. In this paper, we analyze the effects of erroneous load comparisons on the performance of the Po2 scheme. Specifically, we consider {\em load-dependent} and {\em load-independent} errors.
In the load-dependent error model, an incoming job is sent to the server with the larger queue length among the two sampled servers with an error probability $\epsilon$ if the difference in the queue lengths of the two sampled servers is less than or equal to a constant $g$; no error is made if the queue-length difference is higher than $g$. 
For this type of errors, we show that, in the large system limit, the benefits of the Po2 scheme is retained for all values of $g$ and $\epsilon$ as long as the system is heavily loaded, i.e., $\lambda$ is close to $1$. 
In the load-independent error model, 
the incoming job is sent to the sampled server with the {\em maximum load} with an error probability of $\epsilon$ independent of the loads of the sampled servers.
For this model, we show that the performance benefits of the Po2 scheme are retained only if $\epsilon \leq 1/2$; for $\epsilon > 1/2$ we show that the stability region of the system reduces and the system performs poorly in comparison to the {\em random scheme}. 
To prove our stability results, we develop a generic approach to bound the drifts of Lyapunov functions for any state-dependent load balancing scheme. Furthermore, the mean-field analysis in our paper uses a new approach to characterise fixed points which do not admit a recursion. 
\end{abstract}
\section{Introduction}
\label{sec:Intro}

Modern data centres consist of large numbers of parallel servers. Balancing the load across these servers is crucial to ensure better resource utilization and satisfactory quality of service. To distribute the load across the servers uniformly, each incoming job needs to be assigned to an appropriate server in the system.  Assigning an incoming job to one of the servers is done by a job dispatcher using a load balancing scheme. For homogeneous systems, it is well known that the \textit{Join-the-Shortest-Queue} (JSQ) load balancing scheme~\cite{Weber1978,Winston1977optimality}, where an incoming job is assigned to the server having smallest number of ongoing jobs, is optimal in terms of minimizing the average response time of jobs. 
For heterogeneous systems, variants of the JSQ scheme are known to be asymptotically optimal~\cite{sanidhay_performance,sanidhay_wiopt}.
However, implementing JSQ-type schemes is often difficult in practice as dispatchers only have local views of the system and therefore can access the load information of only a subset of servers~\cite{goren2022distributed,goren2021stochastic}. A popular alternative to the JSQ scheme, called the {\em Join-the-Idle-Queue} scheme~\cite{Lu_JIQ_2011,stolyar_JIQ}, requires only the knowledge of the idle servers. However, even this scheme suffers from similar implementation issues as the dispatcher needs to store the idle tokens of a large number of servers.

Due to the challenges mentioned above, randomised schemes, such as the celebrated \textit{Power-of-$d$ choices} or the {\em Pod} scheme~\cite{Mitzenmacherthesis,Vvedenskaya1996}, are widely used in practice~\cite{tarreau2019test,garrett2018nginx}. In the Pod scheme, an arrival is sent to the server with the shortest queue length among a set of $d$ servers sampled uniformly at random. 
For $d=1$, the Pod scheme reduces to the {\em random scheme}. For $d=2$, it is well-known that an exponential improvement in the average response time is achieved in comparison to the random scheme. 
However, this result requires that the server with the least queue length among the two sampled servers is always correctly identified. In a real system, the queue lengths stored at a dispatcher may get outdated due to infrequent updates from the servers~\cite{mike2018netflix,zhou2021asymptotically}. This could result in misidentifying 
the server currently having the smaller queue length.
Another scenario where an error can occur, is when an adversary tries to carry out an attack by misreporting the queue-lengths sent from the servers to the dispatcher~\cite{burke2021misreporting}.
The attacker can manipulate the queue lengths in a way that the dispatcher assigns the job to the server with the maximum load among the two sampled servers. Such erroneous assignments can significantly increase average response time of jobs and may even cause the system to become unstable. The importance of studying the effect of inaccurate load comparisons on load balancing was highlighted as early as 2001 in a survey paper by Mitzenmacher, Richa and Sitaraman~\cite{sitaraman2001power}. However, except from the static setting~\cite{los2022balanced,nadiradze2021achieving}, where there is a finite pool of jobs, this problem has not been studied.
 
{\bf Our Contributions}: This motivates us to consider the effects
that comparison errors can have on the performance of the Po2 scheme in the dynamic setting. 
We consider two types of comparison errors, referred to as the {\em load-dependent} and {\em load-independent} errors.
In the load-dependent error model, an ``error'' is made with probability $\epsilon \in [0,1]$ if the difference in the queue lengths of the sampled servers is sufficiently small (less or equal to a constant $g\geq 0$); if this is not the case (i.e., if the queue-length difference is higher than $g$), then no error occurs. 
An error, in this context, refers to the event where the job is sent to the server having the larger queue length among the two sampled servers. Hence, in this model, the dispatcher makes an error only when the current queue lengths of the sampled servers are close to each other; this is natural to expect when errors occur primarily due to outdated queue lengths at the dispatcher as servers having close queue lengths are likely to be more affected by this type of errors.
To model errors due to adversarial attacks, we consider the load-independent error model in which an error  occurs with probability $\epsilon$ {\em independent} of the current loads of the sampled servers. Clearly, this model of error can have a more drastic impact on the system's performance than the load-dependent error model. 
Our goal is to characterise the performance of the Po2 scheme under these two error models for a system
where there are $n$ unit-rate servers and jobs with exponentially distributed sizes arrive according to a Poisson process with rate $n\lambda$ ($\lambda < 1$). 

It is natural to expect that as $g$ and $\epsilon$ increases (i.e., as the error rate becomes higher), the performance of the Po2 scheme under the load-dependent error model would deteriorate, eventually resulting in a performance poorer than the random scheme. While this is true for light traffic (small values of $\lambda$), we show that, in the heavy traffic limit ($\lambda \to 1$) and large system sizes, the performance of the system remains exponentially better than that under the random scheme for all values of $g$ and $\epsilon$. This implies that the benefits of sampling one additional server in the Po2 scheme outweighs the negative impact of comparison errors when the system operates at its maximum capacity. This result can be explained though the dynamics of the system in the mean-field regime. Specifically, we show the fixed point of the mean-field has a super-exponential decay of tail probabilities for all values of $g$ and $\epsilon$. While this decay rate is slower than that under the standard Po2 scheme, it is still super-exponential and therefore its benefits in comparison to the exponential decay rate under the random scheme become more prominent in the heavy traffic regime.

For the load-independent error model, we show that the benefits of the Po2 scheme are retained only if the error probability $\epsilon$ satisfies $\epsilon \leq 1/2$. For $\epsilon > 1/2$, we show that system becomes unstable for arrival rates larger than $1/2\epsilon$ and the performance becomes worse than that under the random scheme. This can be intuitively explained by the fact that for $\epsilon > 1/2$ the Po2 scheme chooses the server with the larger queue length more often than the server with the smaller queue length. 

From a technical point of view, we make a number of important contributions.
First, we derive the stability region of the system for both error models and establish uniform (in the system size $n$) bounds on the stationary expected queue length per server. These bounds are essential to establish tightness of the stationary measures and interchange of  limits in the mean-field regime. 
The existing results of~\cite{bramson_stability} on JSQ-type load balancing schemes do not apply to our schemes since a job is not always sent to the server with the minimum queue length among the sampled servers. Although the fluid limit results of~\cite{foss_stability} can be applied to derive stability conditions, they do not yield the uniform bounds essential to establish tightness of the stationary distributions. To obtain both stability conditions and the uniform bounds, we use drifts of suitable Lyapunov functions.  However, bounding the drifts of these Lyapunov functions is difficult for our schemes as the schemes compare only a subset of queues at each arrival and do not always choose the shortest queue as the final destination. We develop a generic approach through which the drift can be bounded for any scheme where queue lengths of multiple servers are compared to dispatch each job.

The second important technical contribution is the mean-field analysis of the Po2 scheme under the load-dependent error model. This analysis differs significantly from conventional analysis in that 
the fixed point of the mean-field in this case does not satisfy any recursion. For such a system, even the existence of the  fixed point is not evident. Proving global stability is also more complicated as it relies on induction on the component index. To prove the desired results, we use a new approach that employs bounds on the decay rate of the mean-field process and its monotonicity. We believe that this approach is more broadly applicable to other models where a fixed point does not admit a recursive relationship.

\subsection{Related Works}

In the last two decades, the Pod scheme has emerged as a widely used load balancing scheme due to its promising gains and minimal overhead. It has been studied extensively under various scaling limits and traffic conditions.  
The mean-field scaling limit for this scheme was first studied  for exponential service time distributions in the seminal works~\cite{Mitzenmacherthesis,Vvedenskaya1996}. Their results were later generalised to general service time distributions in~\cite{bramson2012asymptotic,bramson2010randomized}.
The heavy traffic optimality of the Pod scheme has been established in~\cite{chen2012asymptotic,maguluri2014heavy}. In~\cite{mukherjee2020asymptotic}, the analysis of the Pod scheme has been carried out for the case where the number of choices, $d$, is allowed to depend on the system size $n$ and $d(n) =\omega(1)$. In this work, both the mean-field and Halfin-Whitt regimes are considered. In the mean-field limit, the Pod scheme has been shown to reach the same performance as the JSQ scheme. Recently, the Pod scheme has been analysed for different graph topologies. For example, in~\cite{budhiraja2019supermarket}, the Pod scheme is analyzed for non-bipartite graphs and sufficient conditions on the graph sequence is obtained to match the result on complete graphs in the mean-field limit. For the bipartite graphs, the Pod scheme and its variants have been analysed in~\cite{debankur_constrained_2021,zhao2022exploiting}. In all cases, results similar to the complete graph setting have been obtained.
For heterogeneous systems, the Pod scheme has been  studied in~\cite{arpan_tcns,arpan_ssy} where speed-aware versions of the Pod scheme have been shown to yield similar performance benefits.

The above mentioned results for the Pod scheme  
assume that on each arrival the dispatcher has the accurate knowledge of the queue lengths of the $d$ sampled servers. However, this assumption may not be true in practice due to the issues discussed in the introduction. Recently, in~\cite{los2022balanced}, the {\em balls and bins problem} was studied under various noisy load comparison models.
Here, $n$  balls are placed into $n$ bins
sequentially and at each step a ball is placed into a bin from a subset of $d$ bins chosen at random. For the load-dependent error model discussed above, it has been shown that gap between the maximum and the average load is $O(\frac{g}{\log(g)}\log \log(n))$. This result motivates us to consider the effects of noisy load comparisons on the performance of the Po2 scheme in the dynamic setting where jobs are allowed to leave the system after being served.
To the best of our knowledge, this is the first work that studies the Po2 scheme under the erroneous load comparison model in the dynamic setting.


\section{System Model}
\label{sec:system_model}

We consider a system consisting of $n$ parallel servers, 
each with its own queue of infinite buffer size.
Each server is able to process jobs at unit rate.
Jobs  arrive according to a Poisson process with a rate $n\lambda$ ($\lambda < 1$). 
Each job requires a random amount of work, independent and 
exponentially distributed with unit mean. The inter-arrival and job lengths are assumed to be independent of each other. A job dispatcher assigns each incoming job to a queue of a server where jobs are served according to the First-Come-First-Server (FCFS) scheduling discipline. 

The job dispatcher uses the Po2 scheme to assign jobs to the servers\footnote{We use Po2 scheme instead of the Po$d$ scheme with $d\geq 3$ since the gain for $d\geq 3$ is marginal with respect to that for $d=2$.}. Under the classical Po2 scheme, a job is sent to the server with the minimum queue length among two servers, chosen uniformly at random.  However, in practice, the server with the smaller queue length may not be always be correctly identified either due to outdated queue-length information at the dispatcher or due to an attacker  misreporting the queue lengths sent from the servers to the dispatcher.
Motivated by these scenarios,
in this paper, we consider the following versions of the
Po2 scheme where load comparisons are not always accurate. In the following, an error refers to the event where an arrival is sent to the server with the larger queue length among the two sampled servers.


\subsection{Load-Dependent Error Model}

In this model, an error occurs with probability $\epsilon \in [0,1]$ only when the difference in queue lengths of the sampled servers is in the range $(0,g]$ for some constant $g\geq 0$. 
If the queue-length difference is strictly above $g$, then we assume that no error is made, i.e., the job is sent to the server with the smaller queue length. 
In case of a tie, we assume that an arbitrary tie breaking rule based on server indices is used. Without loss of generality (WLOG), we assume that servers are indexed from the index set $[n]=\cbrac{1,2,\ldots,n}$, and, in case of a tie,
the job is sent to the server with the smaller index among the two sampled servers.
We refer to the Po2 scheme under this model of error as the Po2-$(g,\epsilon)$ scheme.

\subsection{Load-Independent Error Model}
  
In this model, an error occurs with probability $\epsilon\in [0,1]$ independent of the current queue lengths of the sampled servers. More precisely, the incoming job is sent to the server having the higher queue length among the sampled servers with probability $\epsilon$ and with probability $1-\epsilon$ it is sent to the server with the smaller queue length. 
Ties are broken in the same way as discussed before.
For simplicity, we refer to the Po2 scheme under this model of error as the Po2-$\epsilon$ scheme. 



\subsection{System State and Notations}

To analyze the system under the schemes discussed above, we first introduce Markovian state descriptors of the system.
We use two Markovian state descriptors.
First, we define the queue-length vector at time $t\geq 0$ as
$$\mf Q^{(n)}(t)=(Q_{k}^{(n)}(t),k \in [n]),$$ 
where $Q_{k}^{(n)}(t)$ denotes
the queue length of the $k^{\textrm{th}}$ server. Second, we define the 
tail measure on the queue lengths at time $t$ as 
$$\mathbf{x}^{(n)}(t)=({x}^{(n)}_{i}(t),i \geq 1),$$
where ${x}^{(n)}_{i}(t)$ denotes the fraction of servers with at least 
$i$ jobs at time $t$. For completeness, we set
$x^{(n)}_{i}(t)=1$ for all $i \leq 0$ and all $t\geq 0$.
From the Poisson arrival and the exponential job size assumption it is clear that both 
$\mf Q^{(n)}=(\mf Q^{(n)}(t), t\geq 0)$ and $\mathbf{x}^{(n)}=(\mathbf{x}^{(n)}(t), t\geq0)$ are Markov processes. When the system is stable, we denote by $\pi_n$ the unique invariant measure of the process $\mf x^{(n)}$ and we use $\mf x^{(n)}(\infty)$
and $\mf Q^{(n)}(\infty)$ to denote the steady-state values of the processes $\mf x^{(n)}$ and $\mf Q^{(n)}$, respectively. As the load balancing scheme does not distinguish between servers, we have 
\begin{equation*}
 \mathbb P(Q^{(n)}_i(t)\geq k)=\E[\indic{ Q^{(n)}_i(t)\geq k}]=(1/n)\sum_{i \in [n]}\E[\indic{Q^{(n)}_i(t)\geq k}]=\E[x^{(n)}_k], \end{equation*}
for each $i \in [n]$ and each $t \in [0,\infty]$.
To define the state space of the process $\mf x^{(n)}$, we first define the space 
$\bar{S}$ as $$\bar{S}\define \{ \mf s=(s_{i}):s_0=1,1\geq s_{i} \geq s_{i+1} \geq 0, \forall  i\geq 1\}.$$
%
%
Note that the space $\bar{S}$ is compact under the norm defined as 
\begin{equation*}
 \norm{ \mf s}= \sup_{i\in \Z_+}\frac{|s_i|}{i+1}, \mf s\in \bar{S}.   
\end{equation*}
%
%
The process $\mf Q^{(n)}$ takes values in $\mb{Z}_+^n$ and
the process $\mf x^{(n)}$ takes values in the space ${S}^{(n)}$
defined as 
$${S}^{(n)}\define \{ \mf s \in \bar{S}: n s_{i}\in \Z_+ \ \forall  i\geq 1\}.$$
%
We further define the space $S$ 
as follows 
$$S\define \{ \mf s\in \bar{S}:\norm{\mf s}_1<\infty\},$$
%
%
where the $\ell_1$-norm, denoted by $\norm{\cdot}_1$, is defined as $\norm{\mf s}_1\define \sum_{i\geq 1} |s_{i}|$
for any $\mf s \in S$. 



\section{Main Results and Insights}
\label{sec:main_results}

In this section, we summarise our main results and discuss their consequences. 
In the following theorem, we characterise the stability region for each of the two schemes discussed above.

\begin{thm}[Stability]
\label{thm:stability}
\begin{temp}
\item For any $g \geq 0, \epsilon \in [0,1]$, and $n \geq 2$ the system
under the Po2-$(g,\epsilon)$ scheme is stable (i.e., the process $\mf x^{(n)}$ is positive recurrent) if and only if $\lambda<1$. 
Furthermore, for $\lambda < 1$, the steady-state average queue length per server is bounded above as 
\begin{equation}
\label{eqn:uniform_bound_g_bounded_1}
\E_{\pi_n}[Q^{(n)}_i(\infty)]=\E_{\pi_n} \sbrac{\sum_{i\geq1}x^{n}_i(\infty)}\leq \frac{(1+g\indic{\epsilon > 1/2})\lambda}{1-\lambda}.
\end{equation}
\label{stability_g_bounded}

\item For any $\epsilon\in[0,1]$ and $n \geq 2$, the system under the Po2-$\epsilon$ scheme is stable, if and only if $\lambda < \min\brac{1,\frac{1}{2\epsilon}}$. 
Furthermore, for $\lambda < \min(1, 1/2\epsilon)$, the steady state average queue length per server is bounded above as 
\begin{equation}
\label{eqn:uniform_bound_epsilon}
\E_{\pi_n}[Q^{(n)}_i(\infty)]=\E_{\pi_n} \sbrac{\sum_{i\geq1}x^{n}_i(\infty)}\leq \frac{\lambda}{1-\max(1,2\epsilon)\lambda}.
\end{equation}
\label{epsilon_stability}
\end{temp}
\end{thm}

In addition to the stability regions, the above theorem gives uniform (in the system size $n$) bounds  on the steady-state mean queue length per server for each scheme. 
These uniform bounds are crucial in establishing the tightness of stationary measures and justifying interchange of the limits in $\lim_{n\to\infty}\lim_{t\to\infty} \mf x^{(n)}(t)=\lim_{t\to\infty}\lim_{n\to\infty} \mf x^{(n)}(t)$, which shows that the mean-field approximation of the steady-state behaviour of the finite system is asymptotically exact.

The bounds in~\eqref{eqn:uniform_bound_g_bounded_1} and~\eqref{eqn:uniform_bound_epsilon} also help us to compare the performance of the Po2-$(g,\epsilon)$ and the Po2-$\epsilon$ schemes to that of the random scheme. For example, when $\epsilon \leq 1/2$, both the upper bounds reduce to $\lambda/(1-\lambda)$ which is the steady-state average queue length per server under the random scheme. This implies that, under both models of error, the Po2 scheme performs better than the random scheme when the error probability $\epsilon \leq 1/2$. This is intuitive, as, for $\epsilon \leq 1/2$, an incoming job under the Po2 scheme is sent to the server with the smaller queue length more often than to the server with the larger queue length. For $\epsilon > 1/2$, however, the schemes may perform poorly in comparison to the random scheme (as both bounds become higher than $\lambda/(1-\lambda)$). 

\begin{figure}[h!]
  \centering
  \includegraphics[width=10cm]{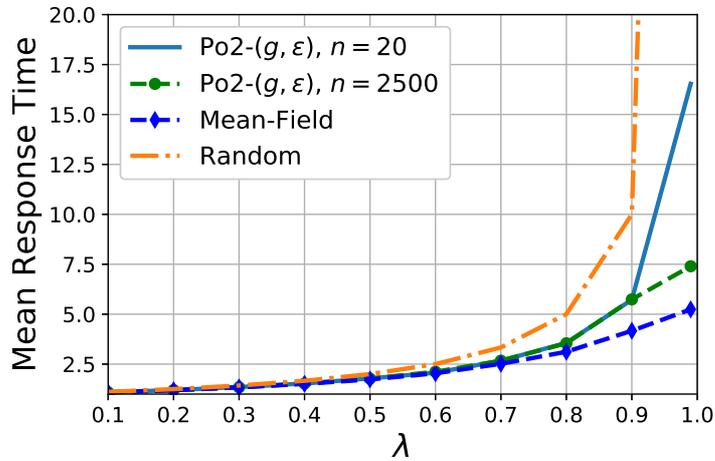}
\caption{Mean response time of jobs under the Po2-$(g,\epsilon)$ scheme as a function of arrival rate $\lambda$ for $\epsilon=0.4$, $g=100$.}
\label{fig:1}
\end{figure}

\begin{figure}[h!]
  \centering
  \includegraphics[width=10cm]{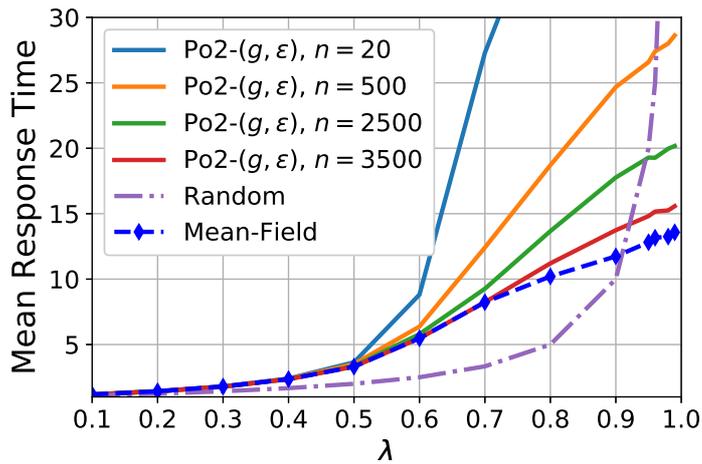}
\caption{Mean response time of jobs under the Po2-$(g,\epsilon)$ scheme and the random scheme as a function of arrival rate $\lambda$ for 
$\epsilon=0.8,g=100$.}
\label{fig:4}
\end{figure}

\begin{figure}[h!]
  \centering
  \includegraphics[width=10cm]{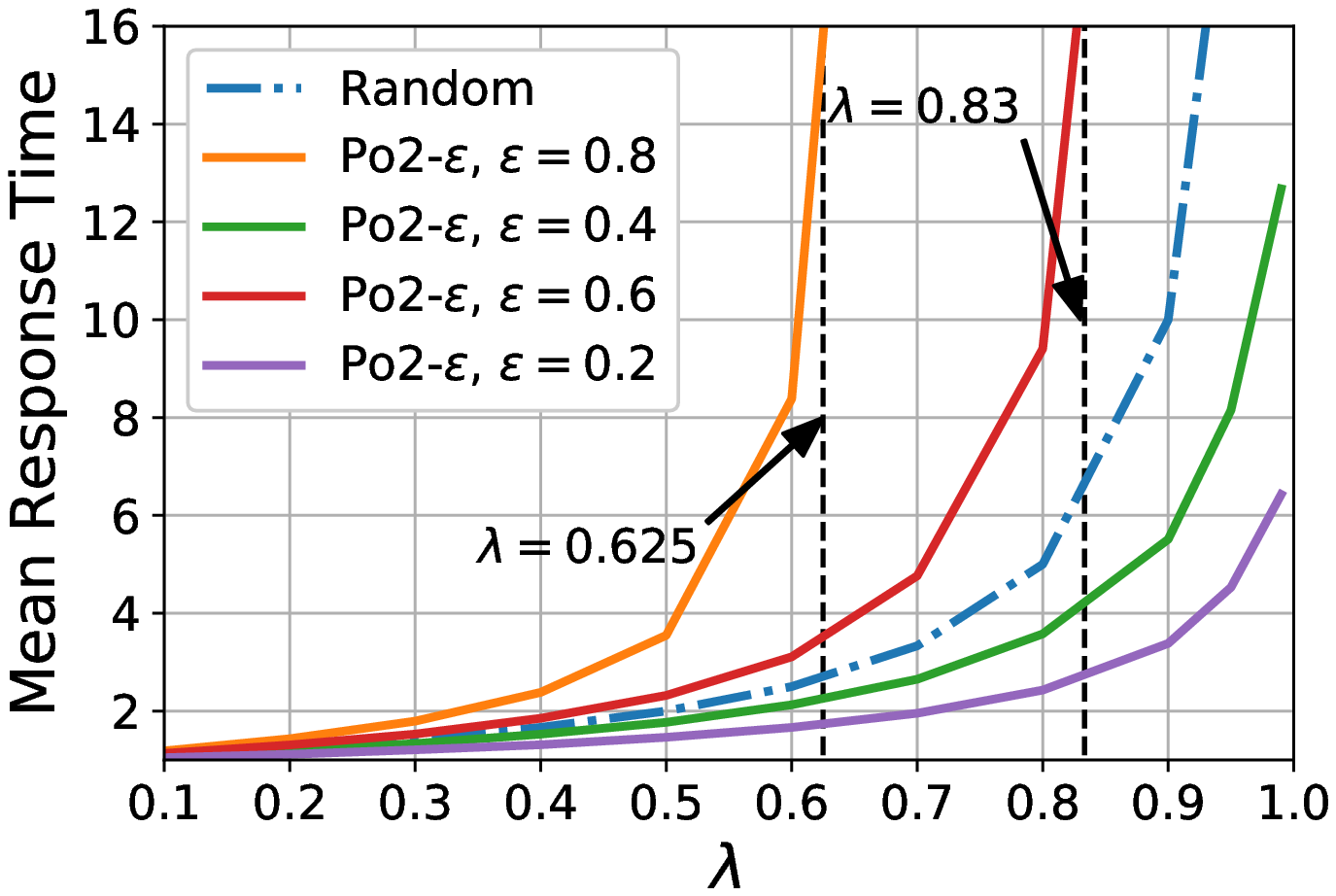}
\caption{Comparison of the Po2-$\epsilon$ scheme for $\epsilon \in \cbrac{0.2,0.4,0.6,0.8}$ with the random scheme. We set $n=200$.}
\label{fig:3}
\end{figure}

This is numerically verified in Figure~\ref{fig:1} and Figure~\ref{fig:4} for the Po2-$(g,\epsilon)$ scheme and in Figure~\ref{fig:3} for the Po2-$\epsilon$ scheme. In each of these figures, we plot the steady-state mean response time of jobs as a function of the normalized arrival rate $\lambda$.
From Figures~\ref{fig:1} and~\ref{fig:3}, we observe that both the schemes outperform the random scheme when $\epsilon \leq 1/2$. For $\epsilon > 1/2$, however, the Po2-$\epsilon$ scheme becomes unstable for $\lambda \geq 1/2\epsilon$ and its performance becomes poorer than that of the random scheme for all $\lambda < 1$. For the Po2-$(g,\epsilon)$ scheme, we observe from Figure~\ref{fig:4} that the system is stable for all $\lambda <1$ even when $\epsilon > 1/2$. However, in this case, the performance of the Po2-$(g,\epsilon)$ scheme is poorer than that of the random scheme for small values of $\lambda$.

The usual approach of proving results similar to the ones stated in Theorem~\ref{thm:stability} consists of coupling and stochastic comparison with the random scheme. However, this approach does not work here since the random scheme can outperform  each of the two schemes when $\epsilon > 1/2$.
Instead, we use drifts of suitable Lyapunov functions to prove Theorem~\ref{thm:stability} which holds for all $\epsilon \in [0,1]$.
Establishing bounds on the drifts of Lyapunov functions is difficult for our schemes as only a subset of servers is compared at each arrival instant and the job is not always sent to the sampled server with the minimum queue length. We develop a generic approach through which the required bounds can be obtained for any scheme where queue lengths of multiple servers are compared to dispatch the incoming jobs.

The asymptotic (in $n$) results for the Po2-$(g,\epsilon)$ scheme are described in the following theorem.

\begin{thm}\textbf{Po2-$(g,\epsilon)$}

\label{thm:main_g_bounded}
\begin{temp}
\item (Mean-Field Limit): Let $g\geq 0$, $\epsilon \in[0,1]$ and assume ${\mathbf{x}}^{(n)}(0)\in  S^{(n)}$ for each $n$ and $ {\mathbf{x}}^{(n)}(0) \overset{a.s}{\to}  \mf u\in S$ under $\ell_1$ as $n \rightarrow \infty$. Then, for each $T \geq 0$, we have $$\sup_{t \in [0,T]}\norm{\mf x^n(t)-\mf x(t)}_1 \overset{a.s}{\to} 0$$ where $\mf x=(\mf x(t)=(x_{i}(t), i \geq 1 ), t \geq 0)$ 
satisfies $\mf x(0)=\mf u$ and for $t\geq0$ and $i\geq1$
\EQN{
\label{eqn:g_bounded_fluid_process}
\dot {x}_{i}(t)=G_i(\mf x(t))\define\lambda  p_{i-1}(\mathbf{x}(t)) - (x_{i}(t)-x_{i+1}(t)).
}
Here for each $i \geq 1$, $G_i$ is the $i^{\textrm{th}}$ component of the function $\mf G=(G_i, i\geq 1): S \to \mb{R}^{\infty}$ and, for $\mf s \in S$, $p_{i-1}(\mf s)$ is defined as
\begin{equation}
\label{eqn:arrival_prop_g_bounded}
p_{i-1}(\mathbf{s})= (s_{i-1}-s_i)\sbrac{ 2\epsilon(s_{i+g} +s_{i-1-g}) + (1-2\epsilon)(s_i+s_{i-1})}.
\end{equation} \label{po2_g_bounded_Process_level_convergence}

\item (Mean-Field Steady State Behaviour): 
For $g\geq 0, \epsilon \in [0,1]$, and $\lambda<1$,
there exists a unique $\mf x^{*} \in S$ such that, $\mf G(\mf x^*)=\mf 0$. In addition, $\mf x^{*}$ satisfies  
\begin{align}
\label{eqn:fx_gb_1}
x_1^{*}&=\lambda,\\
 x_k^{*}&= \lambda  \Big[ 2\epsilon\Big( x_{k-1}^{*} x_{k-1-g}^{*} - \sum_{i=k}^{k+g-1}x_i^{*}(x_{i-1-g}^{*}-x_{i-g}^{*})  \Big)+(1-2\epsilon)(x_{k-1}^{*})^2  \Big], \text{ for } k\geq 2.\label{eqn:fx_gb_2}
\end{align}
Moreover, any solution $\mf x(t)$ of~\eqref{eqn:g_bounded_fluid_process} with $\mf x(0)\in S$ converges to $\mf x^{*}$ in $\ell_1$ as $t \to \infty$.
Therefore, the sequence $(\pi_n)_{n\geq1}$ converges weakly to the Dirac measure $\delta_{\mf x^{*}}$ concentrated on $\mf x^{*}$ as $n \to \infty$.
\label{g_bounded_fixed_point}

\item (Heavy-Traffic Limit): Furthermore, in the heavy traffic as $\lambda \to 1$, we have
\begin{equation}
\label{eqn:perf_gap_g_bounded}
{\limsup}_{\lambda \to 1^-} \frac{T_2^{g,\epsilon}(\lambda)}{\log(T_1(\lambda))}\leq \frac{g+1}{\log(2)},
\end{equation}
where $T_2^{g,\epsilon}(\lambda)=(1/\lambda)\sum_{i \geq 1}x_i^*$
is the steady state limiting (as $n \to \infty$) average response time of jobs under the Po2-$(g,\epsilon)$ scheme and $T_1(\lambda)=1/(1-\lambda)$ is the steady state average response time of jobs under the random scheme.
\label{heavy_gb}


\end{temp}
\end{thm}
%


Hence, the process $\mf x$, defined in Theorem~\ref{thm:main_g_bounded}.(\ref{po2_g_bounded_Process_level_convergence}), characterises the dynamics of the system in the limit as $n \to \infty$. This will be referred to as the {\em mean-field limit} of the system or the {\em mean-field process}.
The evolution of $\mf x$, described in~\eqref{eqn:g_bounded_fluid_process}, can be explained as follows.
For the $n^{\textrm{th}}$ system, the component $x_i^n(t)$ increases by $1/n$ when a job joins a server with queue length exactly $i-1$.  The rate at which this happens is $n\lambda p_{i-1}(\mf x)$, where $p_{i-1}(\mathbf{s})$, for $\mf s \in S$, is the probability that an arrival joins a server with queue-length $i-1$ when the system is in state $\mf s$. The expression of $p_{i-1}(\mf s)$ in~\eqref{eqn:arrival_prop_g_bounded} can be obtained as follows.
Under the Po2-$(g,\epsilon)$ scheme, a job joins a server with queue length $i-1$ under the following scenarios: (1)  One of the two sampled servers is of queue length exactly $i-1$ and the other sampled server is of queue length is at least $i+g$. 
This occurs with probability $2(s_{i-1}-s_i)s_{i+g}$ and, in this case, the job joins the server queue length $i-1$ with probability $1$. (2) One of the sampled server is of queue length $i-1$ and the other sampled server's queue length lies in the range $\cbrac{i-1-g, \dots,i-2}$. This occurs with probability 
$2(s_{i-1} -s_i)(s_{i-1-g}-s_{i-1})$, and, in this case the server with queue length $i-1$ is selected with probability $\epsilon$. (3) One of the sampled server is of queue length $i-1$ and the other sampled server's queue length lies in the range $\cbrac{i+1, \dots,i-1+g}$. This occurs with probability $2(s_{i-1} -s_i)(s_{i}-s_{i+g})$ and, in this case, the server with queue length $i-1$ is selected with probability $(1-\epsilon)$. 4) Finally, both the sampled servers can have the same queue length $i-1$ with probability $(s_{i-1}-s_i)^2$, and, in this case, the job joins a server with queue length $i-1$ with probability $1$. Combining the above probabilities, we obtain the expression for $p_{i-1}(s)$. Similarly, the component $x_i^n(t)$ decreases by $1/n$ when a job leaves a server with queue length $i$ and this occurs with rate $n(x_i-x_{i+1})$. Hence,
the total expected rate of change (drift) in the component $x_i^n(t)$ is given by $G_i(\mf x^n(t))=p_{i-1}(\mathbf{x}^{(n)}(t)) - (x_{i}^{(n)}(t)-x_{i+1}^{(n)}(t))$.
In the limit as $n \to \infty$, this becomes
the rate of change of $x_i(t)$.

In part (\ref{g_bounded_fixed_point}) of Theorem~\ref{thm:main_g_bounded}, we show that, as $t \to \infty$, the mean-field process $\mf x$ converges in $\ell_1$ to the unique point $\mf x^{*} \in S$ 
at which $\mf G(\mf x^*)=\mf 0$.
This point $\mf x^*$ is referred to as the {\em fixed point} of the mean-field since starting
at this point the mean-field remains at this point at all times. Since by part~(\ref{po2_g_bounded_Process_level_convergence}) we have $\mf x^n(t) \to \mf x(t)$ almost surely for each $t \geq 0$, the convergence to the fixed point implies $\lim_{t\to \infty} \lim_{n \to \infty} \mf x^n(t)=\lim_{n\to \infty} \lim_{t \to \infty} \mf x^n(t)=\mf x^*$, which, in turn, means that the fixed point $\mf x^*$ characterises the steady-state behaviour of the system in the limit as $n \to \infty$. In particular, $\lim_{n \to \infty}\mb P(Q^{(n)}_i(\infty))=x^*_i$.

In Theorem~\ref{thm:main_g_bounded}.(\ref{heavy_gb}) we compare the steady-state mean response time of jobs under the Po2-$(g,\epsilon)$ scheme to that under the random scheme when the traffic is high (i.e., $\lambda \to 1$).
Note that~\eqref{eqn:perf_gap_g_bounded} implies that $T_2^{g,\epsilon}(\lambda)=O(\log T_1(\lambda))$ as $\lambda \to 1$. Furthermore, by the previous part of the theorem, the steady-state mean response time of the jobs under the Po2-$(g,\epsilon)$ scheme converges as $n \to \infty$ to $T_2^{g,\epsilon}(\lambda)$. Hence, this result shows that, when the system is heavily loaded, the mean response time of jobs under the Po2-$(g,\epsilon)$ scheme is exponentially smaller than that under the random scheme. This is also verified in Figure~\ref{fig:4} for $\epsilon=0.8$ and $g=100$. Note that for such high error rates, the mean response time of jobs under the Po2-$(g,\epsilon)$ policy can be larger than that underthe random scheme for low values of $\lambda$. However, when $\lambda$ is close to its maximum value $1$, the Po2-$(g,\epsilon)$ scheme performs exponentially better than the random scheme for all values of $g$ and $\epsilon$. This implies that the advantage of having an additional choice in the Po2 scheme outweighs the negative impact the comparison errors when the traffic is high.

The main difficulty in proving Theorem~\ref{thm:main_g_bounded} is that the fixed point $\mf x^*$ cannot be found in closed form.
This is because each component $x_k^*$ in~\eqref{eqn:fx_gb_2} depends not only on the previous components but also on the next $g$ components. This makes it hard to characterise the fixed point; indeed, even the existence of such $\mf x^*$ in $S$ is not evident. This also makes proving the global stability difficult as it uses induction on the component index $k$. To overcome these difficulties, we use the monotonicity of the mean-field and uniform bounds on its tails. We believe that this new approach is generally applicable to similar systems where the fixed point cannot be found in closed form.

We now present the asymptotic results for the Po2-$\epsilon$ scheme in the following theorem.

\begin{thm}\textbf{Po2-$\epsilon$}
\label{thm:Main_thm_epsilon}
\begin{temp}

\item 
(Mean-Field Limit): Let $\epsilon \in[0,1]$ and assume ${\mathbf{x}}^{(n)}(0)\in S^{(n)}$ for each $n$ and $ {\mathbf{x}}^{(n)}(0) \overset{a.s}{\to}  \mf u\in S$ under $\ell_1$ as $n \rightarrow \infty$. Then, for each $T \geq 0$, we have $$\sup_{t \in [0,T]}\norm{\mf x^n(t)-\mf x(t)}_1 \overset{a.s}{\to} 0$$ where $\mf x=(\mf x(t)=(x_{i}(t), i \geq 1 ), t \geq 0)$ 
satisfies $\mf x(0)=\mf u$ and for $t\geq0$ and $i\geq1$
\EQN{
\label{eqn:po2_eps_fluid_process}
\dot {x}_{i}(t)=F_i(\mf x(t))\define \lambda p_{i-1}(\mathbf{x}(t)) -(x_{i}(t)-x_{i+1}(t)).
}
Here for each $i \geq 1$, $F_i$ is the $i^{\textrm{th}}$ component of the function $\mf F=(F_i, i\geq 1): S \to \mb{R}^{\infty}$ and, for $\mf s \in S$, $p_{i-1}(\mf s)$ is defined as
\begin{equation}
\label{eqn:arrival_prop_eps}
p_{i-1}(\mathbf{s})= (1-\epsilon)(s_{i-1}^2 -s_i^2) + \epsilon ( (1-s_i)^2 - (1-s_{i-1})^2 ).
\end{equation}
We refer to the process $\mf x$ as the mean-field limit of the sequence $(\mf x^{(n)})_{n\geq1}$.
\label{po2_eps_Process_level_convergence}
\item (Mean-Field Steady State Behaviour): For $\epsilon \in [0,1]$  and $\lambda < \min(\frac{1}{2\epsilon},1)$, there exits $\mf x^*\in S$ such that, if $\mf x(0)=\mf x^*$, then $\mf x(t)=\mf x^*$ for all $t\geq 0$. Furthermore,  $\mf x^*$ satisfies the following recursion
\EQN{
\label{eqn:po2_eps_Fixed_Point}
\begin{aligned}
x_1^*=\lambda, \ \ \ x^*_{i}=\lambda \sbrac{(1-2\epsilon)(x_{i-1}^*)^2 + 2\epsilon x_{i-1}^*} \ \ \forall \ i\geq2.
\end{aligned}
}

%
Moreover, any solution $\mf x(t)$ of~\eqref{eqn:po2_eps_fluid_process} with
$\mf x(0) \in S$ converges to $\mf x^*$ in $\ell_1$, i.e., $\norm{\mf x(t) - \mf x^*}_1 \to 0$ as $t \to \infty$.
The above results imply that the sequence $(\pi_n)_n$ converges weakly to the Dirac measure $\pi_*=\delta_{\mf x^*}$ concentrated on $\mf x^*$ as $n \rightarrow \infty$.
\label{fixed_point}
\item (Heavy-Traffic Limit):
For $\epsilon\leq 1/2$ and $\lambda<1$, we have 
\begin{equation}
\label{eqn:perf_gap}
\lim_{\lambda \to 1^-} \frac{T_2^{\epsilon}(\lambda)}{\log(T_1(\lambda))}=\frac{1}{\log(2-2\epsilon)},
\end{equation} 
where $T_2^{\epsilon}(\lambda)=\frac{\sum_{k=1}^{\infty} x_k^*}{\lambda}$ is the limiting (as $n \to \infty$) steady state average response time of jobs under the Po2-$\epsilon$ scheme.
\label{comp_ratio_epsilon}
\end{temp}
\end{thm}

In parts (\ref{po2_eps_Process_level_convergence}) and
(\ref{fixed_point})
of Theorem~\ref{thm:Main_thm_epsilon}, we characterize the mean-field limit $\mf x$ and its fixed point $\mf x^*$ under the Po2-$\epsilon$ scheme. 
As before, we show that the fixed point is unique and globally asymptotically stable. In the last part (part \ref{comp_ratio_epsilon}) of Theorem~\ref{thm:Main_thm_epsilon}, we compare the mean response time of jobs under the Po2-$\epsilon$ scheme to that under the random scheme in the limit as $n \to \infty$. Our result indicates that when $n$ is large and $\lambda$ is close to $1$, the steady state mean response time of jobs under the Po2-$\epsilon$ scheme satisfies $T_2^{\epsilon}(\lambda) \approx c_{\epsilon}\log(T_1(\lambda)$, where $c_\epsilon=1/(\log(2-2\epsilon))$. This means that an exponential reduction in the steady state mean response time is achieved as long as $\epsilon \leq 1/2$. Hence, the Po2-$\epsilon$ scheme retains the benefits of the Po2 scheme as long as $\epsilon \leq 1/2$.

\section{Stability and Uniform Bounds}
\label{sec:stability}

In this section, we find the stability regions for the Po2-$(g,\epsilon)$ and the Po2-$\epsilon$ schemes and derive uniform bounds on the steady-state queue length per server (Theorem~\ref{thm:stability}) using drifts of appropriate Lyapunov functions. 
We first develop a general framework to analyse any load balancing scheme that compares the queue lengths of two uniformly sampled servers to dispatch every job. Note that it is easy to generalise this framework further to cases where more than two servers are sampled and the sampling is not necessarily uniform.

For any function $V:\mb{Z}_+^n \to [0,\infty)$, the drift of $D_{\mf Q^n} V$ is defined as the expected rate of change in the value of the function along the trajectory of the process $\mf Q^n$ given the current state. More precisely,
\begin{align}
    D_{\mf Q^{(n)}}  V(\mf Q) &= \lim_{h \to 0}\frac{1}{h}\E[V(\mf Q^{(n)}(t+h))-V(\mf Q^{(n)}(t))| \mf Q^n(t)=\mf Q] \nonumber\\
    &= \sum_{i=1}^n [ r_{i}^{+,n}(\mf Q)(V(\mf Q+\mf e_{i}^{(n)})-V(\mf Q))+r_{i}^{-,n}(\mf Q)(V(\mf Q-\mf e_{i}^{(n)})- V(\mf Q))],
     \label{eq:genq_def}
\end{align}
where $\mf e_{i}^{(n)}$ denotes the $n$-dimensional 
unit vector with one
in the $i^{\textrm{th}}$ position; 
$r^{\pm,n}_{i}(\mf Q )$
are the transition rates from the state $\mf Q $ to the states $\mf Q \pm \mf e_{i}^{(n)}$. According to the Foster-Lyapunov theorem (Proposition D.3 of \cite{Kelly_book}), to prove the stability or positive recurrence of the process $\mf Q^{(n)}$, it is sufficient to show the existence of at least one function $V:\mb{Z}_+^n \to [0,\infty)$ such that $V(\mf Q) \to \infty$ when $\norm{\mf Q}_1 \to \infty$ and $D_{\mf Q^n} V(\mf Q) < 0$ for all states $\mf Q$ lying outside a compact subset of the state-space. To further obtain uniform bounds on the stationary queue lengths, we use the fact (Proposition 1 of~\cite{Glynn_bounds}) that $\E_{\pi_n}[D_{\mf Q^{(n)}} V(\mf Q(\infty))] \geq 0$ if $D_{\mf Q^{(n)}} V(\mf Q)$ is uniformly bounded 
for all states $\mf{Q} \in \mb{Z}_+^n$. 

The rate of departure from the $i^{\textrm{th}}$ server is given by $r_{i}^{-,n}(\mf Q)=\indic{Q_i>0}$. 

For any scheme which compares the states of two servers to dispatch the job to one of the servers, we define the class of an arrival as the (unordered) pair $\brac{i,j}$ of servers sampled at the arrival instant. Let $\mathcal C$ denote the collection of all such classes. Since $|\mathcal{C}|=\binom{n}{2}$ and a job is equally likely to belong to one of these classes, the arrival rate of any class $(i,j)\in \mc C$ is $n\lambda/\binom{n}{2}=2\lambda/(n-1)$.
Hence, we can write the rate of arrival to the $i^{\textrm{th}}$ server as
\begin{equation}
\label{eqn:r+}
r_{i}^{+,n}(\mf Q)=\frac{2\lambda}{n-1}\sum_{j\in[n],j\neq i}p(Q_i,Q_j),
\end{equation}
where $p(Q_i,Q_j)$ is the probability that a class $\brac{i,j}$ job joins the server $i$ when the queue lengths of servers $i$ and $j$ are $Q_i$ and $Q_j$, respectively. 
Note that the probability $p(Q_i,Q_j)$ depends on the load balancing scheme used by the dispatcher. The exact expression of $p(Q_i,Q_j)$ for each scheme is given later, but it is important to note that  $p(Q_i,Q_j)+ p(Q_j,Q_i)=1$ since a class $(i,j)$ job joins either server $i$ or server $j$ with probability $1$. 

Now, for the Lyapunov function $V:\mb{Z}_+^n \to [0,\infty)$ defined as 
$$ V(\mf Q)= \sum_{i =1}^n Q_{i}^2,$$ the drift given in~\eqref{eq:genq_def} simplifies to
\begin{align}
D_{\mf Q^{(n)}} V(\mf Q) 
    =  \sum_{i =1}^n \cbrac{2[r_{i}^{+,n}(\mf Q) -r_{k}^{-,n}(\mf Q) ] Q_{i} +   [r_{i}^{+,n}(\mf Q) +r_{i}^{-,n}(\mf Q)]},
    \nonumber
\end{align}
which upon further simplification gives
\begin{equation}
\label{eqn:D_drift}
 D_{\mf Q^{(n)}} V(\mf Q) 
    =  2\sum_{i =1}^n r_{i}^{+,n}(\mf Q)Q_{i} -2\sum_{i\in[n]}Q_{i} + n\lambda +B(\mf Q),
\end{equation}
where $B(\mf Q)=\sum_{i\in [n]}\indic{Q_i >0}$ represents the number of busy servers when system is in state $\mf Q$. In the above, we have used the facts $r_i^{-,n}(\mf Q)=\indic{Q_i > 0}$ and $\sum_{i\in[n]}r_{i}^{+,n}(\mf Q)=n \lambda$.
Moreover, using~\eqref{eqn:r+}, the first term in the exresssion of the drift can be written as  
\begin{align}
\label{eqn:r_i+_sum}
&\sum_{i =1}^n r_{i}^{+,n}(\mf Q)Q_{i}=  \frac{2\lambda}{(n-1)} \sum_{i\in[n]}  \sum_{j\in[n],j\neq i}Q_ip(Q_i,Q_j) \nonumber\\
&=\frac{2\lambda}{(n-1)}  \sum_{(i,j) \in \mc{C}}(Q_ip(Q_i,Q_j)+Q_jp(Q_j,Q_i)).
\end{align}
Thus, to obtain the stability region and the uniform bound on steady-state queue length, we need to obtain upper bounds on $Q_i p(Q_i,Q_j)+Q_j p(Q_j,Q_i)$ for each scheme.

\subsection{Po2-$(g,\epsilon)$ Scheme}

For the Po2-$(g,\epsilon)$ scheme, the probability $p(Q_i,Q_j)$ for any class $(i,j)\in \mc{C}$ is given by
\begin{multline}
\label{eqn:g_bounded_prob}
p(Q_i,Q_j)=\indic{Q_j-Q_i\geq g+1} +(1-\epsilon)\indic{Q_j-Q_i\in (0,g]}\\+\epsilon \indic{Q_i-Q_j \in (0,g]}+\indic{Q_i=Q_j,i<j}.
\end{multline}
Using the above expression, we obtain the following bound for Po2-$(g,\epsilon)$ scheme. 

\begin{lem}
\label{lem:po2_g_ep_bound_l}
For $g\geq0$, $\epsilon\in[0,1]$, and for any class $(i,j)\in \mc{C}$, under the Po2-$(g,\epsilon)$ scheme, we have 
\begin{equation}
Q_ip(Q_i,Q_j)+Q_jp(Q_j,Q_i) \leq \frac{Q_i+Q_j}{2} +g \indic{\epsilon>1/2}.
\end{equation}
\end{lem}
\begin{proof}
To prove the lemma, we first observe that for any $a\leq b$, and $w_1=1-w_2 \in [1/2,1]$ we have 
\begin{equation}
\label{eqn:weight_mean}
w_1 a+w_2b \leq \frac{a+b}{2}.
\end{equation}
The result of the lemma is direct when $Q_i=Q_j$. So, we consider the case $Q_i<Q_j$. Note that the proof for $Q_i > Q_j$ is exactly the same with $Q_i$ and $Q_j$ interchanged. For $\epsilon \leq 1/2$,  using~\eqref{eqn:g_bounded_prob} and $Q_i < Q_j$, we have 
\begin{align*}
 p(Q_i,Q_j)&=\indic{Q_j-Q_i\geq g+1} +(1-\epsilon)\indic{Q_j-Q_i\in (0,g]},\\
 &\geq \frac{1}{2}\indic{Q_j-Q_i\geq g+1} +\frac{1}{2}\indic{Q_j-Q_i\in (0,g]}=\frac{1}{2},
\end{align*}
where the last equality follows because of the assumption $Q_i < Q_j$. Therefore, $p(Q_i,Q_j)\in [1/2,1]$. Since $p(Q_i,Q_j)=1-p(Q_j,Q_i)$, using~\eqref{eqn:weight_mean}, we have
$$Q_ip(Q_i,Q_j)+Q_jp(Q_j,Q_i) \leq (Q_i+Q_j)/{2}.$$
Next, we note from~\eqref{eqn:g_bounded_prob} and $Q_i < Q_j$ that
\begin{align*}
\label{eqn:lem_temp_1_g}
 Q_ip(Q_i,Q_j)&=Q_i \brac{\indic{Q_j-Q_i\geq g+1} +(1-\epsilon)\indic{Q_j-Q_i\in (0,g]}} \nonumber\\
 &=Q_i \brac{(1-\epsilon)+\epsilon \indic{Q_j-Q_i\geq g+1}},
 \end{align*}
which gives
\begin{multline*}
Q_ip(Q_i,Q_j)+Q_jp(Q_j,Q_i) =Q_i \brac{(1-\epsilon)+\epsilon \indic{Q_j-Q_i\geq g+1}}
+Q_j \epsilon \indic{Q_j-Q_i \in (0,g]}.
\end{multline*}
Therefore, using 
$$Q_j \indic{Q_j-Q_i \in (0,g]} \leq (g+ Q_i) \indic{Q_j-Q_i \in (0,g]} \leq (g+ (Q_i+Q_j)/2) \indic{Q_j-Q_i \in (0,g]},$$ and $Q_i < (Q_i+Q_j)/2 +g$ in the above, we obtain 
\begin{equation*}
Q_ip(Q_i,Q_j)+Q_jp(Q_j,Q_i) \leq g+\frac{Q_i+Q_j}{2}.
\end{equation*}
This completes the proof.
\end{proof}

{Proof of Theorem~\ref{thm:stability}.(\ref{stability_g_bounded})}: Using the bound of Lemma~\ref{lem:po2_g_ep_bound_l}, the RHS of~\eqref{eqn:r_i+_sum} can be bounded as
\begin{align}
\label{eqn:r_+_g_bounded_upper_bound}
&\sum_{i =1}^n r_{i}^{+,n}(\mf Q)Q_{i} \leq \frac{2\lambda}{(n-1)} \sum_{(i,j) \in \mc {C}}\cbrac{\frac{Q_i+Q_j}{2} +g \indic{\epsilon>1/2}} \nonumber\\
&=\frac{2\lambda}{(n-1)} \brac{\frac{(n-1)}{2}\sum_{i\in[n]}Q_i +  g\indic{\epsilon>1/2} \binom{n}{2}}\nonumber\\
&=\lambda \sum_{i\in[n]}Q_i +  n\lambda g\indic{\epsilon>1/2} 
\end{align}
Therefore, using~\eqref{eqn:r_+_g_bounded_upper_bound} in ~\eqref{eqn:D_drift}, we can upper bound the drift $ D_{\mf Q^{(n)}} V(\mf Q)$ for the Po2-$(g,\epsilon)$ scheme as 
\begin{equation}
\label{eqn:drift_g_eps_final}
 D_{\mf Q^{(n)}} V(\mf Q) 
    \leq  2(\lambda-1)\sum_{i=1}^nQ_i+2 n\lambda g\indic{\epsilon>\frac{1}{2}}+ n\lambda +B(\mf Q).
\end{equation}
Now, since $B(\mf Q)\leq n$ and $\lambda<1$, the drift $ D_{\mf Q^{(n)}} V(\mf Q) $ is strictly negative whenever $\sum_{i\in [n]}Q_i > n( \lambda(2g \indic{\epsilon>1/2}+1)+1)/2(1-\lambda)$, and is bounded above by $n( \lambda(2g \indic{\epsilon>1/2}+1)+1)$, otherwise. This shows that the system under the Po2-$(g,\epsilon)$ scheme is stable for all $\lambda < 1$. The necessity of this condition for stability can be established easily by showing that the drift of the Lyapunov function $V_1(\mf Q)=\sum_{i\in [n]}Q_i$ is always non-negative when $\lambda \geq 1$.

To prove~\eqref{eqn:uniform_bound_g_bounded_1}, recall from the previous paragraph that $D_{\mf Q^{(n)}} V(\mf Q)$ is uniformly bounded by $n(\lambda(2g \indic{\epsilon>1/2}+1)+1)$ for all $\mf Q\in \mb{Z}_+^n$.
%
%
This implies that
$\E_{\pi_n} \sbrac{D_{\mf Q^{(n)}} V(\mf Q^{(n)}(\infty))}\geq0$.
Therefore, taking expectation of~\eqref{eqn:drift_g_eps_final}
and using the rate conservation equation $\E_{\pi_n}[B(\mf Q)]=n\lambda$ (which holds in steady-state), we obtain
\begin{align*}
\frac{n \lambda (1+g \indic{\epsilon>1/2})}{(1-\lambda)}\geq \E_{\pi_n} \sbrac{\sum_{i\in[n]}Q_i^{(n)}(\infty)}=n\E_{\pi_n} \sbrac{Q_i^{(n)}(\infty)},
\end{align*}
where last equality follows due to the exchangeability of  $\pi_n$. \qed

\subsection{Po2-$\epsilon$ Scheme}

For the Po2-$\epsilon$ scheme, 
$p(Q_i,Q_j)$ for any class $(i,j)\in \mc{C}$ is given by
\begin{equation}
\label{eqn:po2_eps_prob}
p(Q_i,Q_j)=(1-\epsilon)\indic{Q_i<Q_j}+\epsilon \indic{Q_i>Q_j}+\indic{Q_i=Q_j,i<j}.
\end{equation}
Using the expression above, we obtain the following bound.

\begin{lem}
\label{lem:po2_ep_bound_l}
For $\epsilon\in[0,1]$ and for any class $(i,j)\in \mc{C}$, under the Po2-$\epsilon$, scheme we have 
\begin{equation}
Q_ip(Q_i,Q_j)+Q_jp(Q_j,Q_i) \leq (Q_i+Q_j) \max(1/2,\epsilon).
\end{equation}
\end{lem}
\begin{proof}
For $Q_i=Q_j$ the above inequality follows directly.
Similar to the proof of Lemma~\ref{lem:po2_g_ep_bound_l}, it is sufficient to  consider the case $Q_i<Q_j$. For $\epsilon\leq 1/2$, using~\eqref{eqn:po2_eps_prob} we have $p(Q_i,Q_j)=(1-\epsilon)\in [1/2,1]$. Therefore, from~\eqref{eqn:weight_mean}, it follows that $Q_ip(Q_i,Q_j)+Q_jp(Q_j,Q_i) \leq (Q_i+Q_j)/{2}$. When $\epsilon>1/2$, for $Q_i<Q_j$, we can write 
\begin{align*}
Q_ip(Q_i,Q_j)+Q_jp(Q_j,Q_i)=Q_i (1-\epsilon)+Q_j\epsilon \leq \epsilon (Q_i+Q_j).
\end{align*}
Hence, the proof is complete.
\end{proof}

{\em Proof of Theorem~\ref{thm:stability}.(\ref{epsilon_stability})}: Using the bound of Lemma~\ref{lem:po2_ep_bound_l}, the RHS of~\eqref{eqn:r_i+_sum} can be bounded as 
\begin{align}
\label{eqn:r_+_eps_bounded_upper_bound}
&\sum_{i =1}^n r_{i}^{+,n}(\mf Q)Q_{i} \leq \frac{2\lambda}{(n-1)} \sum_{(i,j) \in \mc {C}}(Q_i+Q_j)\max(1/2,\epsilon), \nonumber\\
& =\lambda\max(1,2\epsilon)\sum_{i\in[n]}Q_i.
\end{align}
Therefore, using~\eqref{eqn:r_+_eps_bounded_upper_bound} in~\eqref{eqn:D_drift}, we upper-bound the drift $D_{\mf Q^{(n)}} V(\mf Q)$ for the Po2-$\epsilon$ scheme as
\begin{equation}
\label{eqn:drift_eps_final}
 D_{\mf Q^{(n)}} V(\mf Q) 
    \leq  2(\lambda\max(1,2\epsilon)-1)\sum_{i\in [n]}Q_i+ n\lambda +B(\mf Q).
\end{equation}
Since $B(\mf Q)\leq n$, the above implies that, for $\lambda < \min(1,1/2\epsilon)$, the drift is strictly negative whenever $\sum_{i\in [n]}Q_i > n(\lambda+1)/(1-\lambda \max (1,2\epsilon))$. This shows that the system under Po2-$\epsilon$ scheme is stable for all $\lambda<\min(1,1/2\epsilon)$. Furthermore, since $\sup_{\mf Q \in \Z_+^n} D_{\mf Q^{(n)}} V(\mf Q)\leq n(\lambda+1)$,
using $\E_{\pi_n} \sbrac{D_{\mf Q^{(n)}} V(\mf Q^{(n)}(\infty))}\geq0$, we obtain
\begin{align*}
\frac{n \lambda }{(1-\lambda \max(1,2\epsilon))}\geq \E_{\pi_n} \sbrac{\sum_{i\in[n]}Q_i^{(n)}(\infty)}=n\E_{\pi_n} \sbrac{Q_i^{(n)}(\infty)},
\end{align*}
which proves~\eqref{eqn:uniform_bound_epsilon}.

Next we prove that for $\lambda \geq \min(1,\frac{1}{2\epsilon})$ the system is unstable. For $\epsilon\leq 1/2$ and $\lambda\geq 1$, the process $\mf Q^{(n)}$ is not positive recurrent. This follows using the same argument as used in the stability proof of the Po2-$(g,\epsilon)$ scheme. Now, for $\epsilon>1/2$ and $2\lambda \epsilon>1$, we consider the Lyapunov function 
$$V_2(\mf Q)=  Q_{i^*(\mf Q)},$$
where $i^*(\mf Q)=\arg\max_{i \in [n]}Q_i$ and $i^*(\mf Q)$ is the minimum such index.
Using~\eqref{eq:genq_def},  the drift of the function $V_2(\mf Q)$ can be written as 
\begin{align}
D_{\mf Q^{(n)}} V_2(\mf Q) =  \brac{ r_{i^*(\mf Q)}^{+,n}(\mf Q) -\indic{Q_{i^*(\mf Q)}>0}}\geq \brac{ r_{i^*(\mf Q)}^{+,n}(\mf Q) -1}.\nonumber
\end{align}
From~\eqref{eqn:r+}, we have
\begin{equation}
r_{i^*(\mf Q)}^{+,n}(\mf Q) = \frac{2\lambda}{n-1}\sum_{j \neq i^*(\mf Q)} p(Q_i,Q_j)\geq \frac{2\lambda}{n-1} \sum_{j \neq i^*(\mf Q)} \epsilon=2\epsilon\lambda.\nonumber
\end{equation}
where the inequality follows from~\eqref{eqn:po2_eps_prob} since
$Q_{i^*(\mf Q)}\geq Q_j$ for any $j\in [n]$ and $i^*(\mf Q)$ is the minimum such index. Hence, $D_{\mf Q^{(n)}} V_2(\mf Q)\geq 2\epsilon\lambda-1\geq 0$
for all $\mf Q \in \mb{Z}_+^n$.
Furthermore, since $D_{\mf Q^{(n)}} V_2(\mf Q) \leq n \lambda$, the result follows from the Foster-Lyapunov criterion for transience and null recurrence (Theorem 3.3.10 of \cite{srikant2013communication}). \qed
\section{Mean-Field Analysis of the Po2-$(g,\epsilon)$ scheme}
\label{sec:po2_g_bounded_analysis}

In this section, we prove the main results for the Po2-$(g,\epsilon)$ scheme stated in Theorem~\ref{thm:main_g_bounded}.

\subsection{Mean-Field Limit of the  Po2-$(g,\epsilon)$ Scheme}

First, we establish the mean-field limit of the Po2-$(g,\epsilon)$ scheme (Theorem~\ref{thm:main_g_bounded}.(\ref{po2_g_bounded_Process_level_convergence})) for any $g\geq0$ and $\epsilon \in [0,1]$. Note that under the Po2-$(g,\epsilon)$ scheme, the rate of transition of the process $\mf x^{(n)}$ from state $\mf x \in S^{(n)}$ to state $\mf y \in S^{(n)}$ is given by
\EQN{
\label{eqn:rate_transition_g_bounded}
r^{(n)}_{\mf x, \mf y} =
\begin{cases}
                                      n \lambda p_{i-1}(\mf x),& \text{if $\mf y=\mf x + \mf e_i/n$} \\
                                   n(x_i-x_{i+1}),& \text{if $\mf y= \mf x -\mf e_i/n$}
\end{cases}, \ \forall i\geq1,
}
where $p_{i-1}(\mf x)$ is as defined in~\eqref{eqn:arrival_prop_g_bounded} and
 $\mf e_i$ is the $i^{\textrm{th}}$ unit vector in $\mb R^{\infty}$.
Clearly, the above rates satisfy the transition structure for
a {\em density-dependent jump Markov chain}~\cite{kurtz_book,Mitzenmacherthesis}.
Furthermore, it is easy to verify that $\sum_{y \in S} r^{(n)}_{\mf x, \mf y} < n(\lambda +1)$ for all $x \in S$ and
the function $\mf G: S \to \mb R^{\infty}$
%
is Lipschitz under the $\ell_1$-norm with a Lipschitz constant of $L_{\lambda}^{\epsilon}=\lambda(16 \epsilon+4)+2$ (proved in Lemma~\ref{lem:lip_g_bounded}). Hence, using the Kurtz's theorem for density-dependent jump Markov processes [\cite{kurtz1981approximation}, Chapter 8], we obtain the desired result.
\begin{lem}
\label{lem:lip_g_bounded}
The function $\mf G(\mf x)$ is Lipschitz under the $\ell_1$ norm with constant $L_{\lambda}^{\epsilon}=\lambda(16 \epsilon+4)+2$. 
\end{lem}
\begin{proof}
Note that we can write~\eqref{eqn:arrival_prop_g_bounded} for any $\mf x \in S$ as 
\begin{equation}
\label{eqn:prob_g_bounded_form_2}
p_{i-1}(\mf x)=\epsilon \sbrac{(x_i-x_{i+g})^2+ (x_{i-1-g}-x_i)^2-(x_{i-1-g} -x_{i-1})^2 - (x_{i-1}-x_{i+g})^2} +x_{i-1}^2 -x_i^2.
\end{equation}
Let $\mf x, \mf y \in S$. Then we have  
 \begin{align}
 \label{eqn:gb_lip_1}
 \norm{\mf G(\mf x) -\mf G(\mf y)}_1&=\sum_{i\geq 1} |  \lambda p_{i-1}(\mathbf{x}) - (x_{i}-x_{i+1}) -  \lambda p_{i-1}(\mathbf{y}) + (y_{i}-y_{i+1}) | \nonumber\\
 &\leq \lambda\sum_{i\geq 1}|p_{i-1}(\mathbf{x})-p_{i-1}(\mathbf{y})| + 2 \norm{ \mf x -\mf y}_1.
 \end{align}
 Now from~\eqref{eqn:prob_g_bounded_form_2}, using the triangle inequality we can write for any $\mf x, \mf y \in S$
 \begin{multline*}
|p_{i-1}(\mathbf{x})-p_{i-1}(\mathbf{y})| \leq \epsilon \Big\{ | (x_i-x_{i+g})^2 - (y_i-y_{i+g})^2  | + | (x_{i-1-g}-x_{i})^2 - (y_{i-1-g}-y_{i})^2  |\\
+ | (y_{i-1}-y_{i+g})^2 - (x_{i-1}-x_{i+g})^2  |+ | (y_{i-1-g}-y_{i-1})^2 - (x_{i-1-g}-x_{i-1})^2  |\Big \}\\ + |x_{i-1}^2-y_{i-1}^2|+ |y_i^2-x_i^2|.
 \end{multline*}
 Therefore, using the above inequality we can write
 \begin{align}
\label{eqn:gb_lip_2}
 \lambda\sum_{i\geq 1}|p_{i-1}(\mathbf{x})-p_{i-1}(\mathbf{y})| \leq \lambda(16\epsilon +4) \norm{\mf x -\mf y}_1.
 \end{align}
 Hence, the result follows from~\eqref{eqn:gb_lip_1} and~\eqref{eqn:gb_lip_2}.
\end{proof}
\subsection{Mean-Field Steady State Behaviour for the Po2-$(g,\epsilon)$ Scheme}

We now turn to the proof of Theorem~\ref{thm:main_g_bounded}.(\ref{g_bounded_fixed_point}) which shows that the mean-field process $\mf x$ given by~\eqref{eqn:g_bounded_fluid_process} has a unique fixed point $\mf x^*$ which satisfies~\eqref{eqn:fx_gb_1} and~\eqref{eqn:fx_gb_2}. Moreover, we show that the fixed point $\mf x^*$ is globally stable, i.e., all trajectories of the mean-field process $\mf x$ starting in $S$ converges to $\mf x^*$. 

For any $\mf u \in S$, let $\mf x(t,\mf u)$ denote the trajectory of the mean-field process starting at state $\mf u$. 
Further, define $v_{k}(t,\mf u)=\sum_{i\geq k}x_{i}(t,\mf u)$ and $v_{k}(\mf u)=\sum_{i\geq k}u_{i}$
for each $k\geq1$.
When the context is clear, we shall drop the dependence 
of the trajectory on the initial state $\mf u$ and on the time $t$.


\begin{lem}
\label{lem:g_bounded_ODE_Solution_Monotone}
Let $g\geq 0, \epsilon\in [0,1]$. The following statements hold for the process $\mf x$ defined in Theorem~\ref{thm:main_g_bounded}. 
\begin{enumerate}
\item If $\mathbf{u}\in S$, then, for any $t\geq 0$, we have $\mf x(t,\mf u) \in S$.
\item For any $\mf u, \mf u' \in S$ satisfying $\mf u \leq \mf u'$
we have ${\mf x}(t, \mf u) \leq {\mf x}(t, \mf u')$ for all $t \geq 0$, where the inequality $\leq$ is understood component-wise.
\end{enumerate}
\end{lem}

\begin{proof}
We first note from~\eqref{eqn:g_bounded_fluid_process} and~\eqref{eqn:arrival_prop_g_bounded} that 
\begin{align}
 \dot x_i \leq \lambda p_{i-1}(\mf x) &\leq \lambda [2\epsilon (x_{i-1-g}^2-x_i^2)+(1-2\epsilon)(x^2_{i-1}-x_i^2)]\nonumber\\
 &\leq \lambda \max(1,2\epsilon)x_{i-1-g}^2
\leq M_{\lambda,\epsilon}x_{i-1-g},\nonumber
\end{align}
where $M_{\lambda,\epsilon}=\lambda\max(1,2\epsilon)$. This implies that for each $i \geq 1$ we have $$x_i(t) \leq x_i(0)+M_{\lambda,\epsilon}\int_{0}^t x_{i-1-g}(s) ds.$$
Using the above recursively for each $i \geq 1$, we obtain
\begin{equation}
\label{eqn:g_bounded_sum}
x_i(t) \leq x_i(0)+\sum_{k=0}^{i-g-1}x_{k}(0)\frac{(M_{\lambda,\epsilon}t)^{i-g-k}}{(i-g-k)!}.
\end{equation}
Summing the above for all $i \geq 1$
we obtain $$v_1(\mf x(t))\leq (g+1+v_1(\mf x(0)))\exp(M_{\lambda,\epsilon}t).$$
This shows that if $v_1(\mf x(0)) < \infty$, then $v_1(\mf x(t)) < \infty$ for all $t$, thus establishing the first part of the lemma.

For the second part, we note that  for each $i\geq 1$, $G_i(\mf x)=\lambda p_{i-1}(\mf x)-(x_i-x_{i+1})$ is non-decreasing with respect to all $x_k$, $k\neq i$. Hence the result follows from Theorem 5.3 of~\cite{Deimling2006}. 
\end{proof}

The second property stated in Lemma~\ref{lem:g_bounded_ODE_Solution_Monotone} is called the {\em quasi-monotonicity} of the process $\mf x$. This property ensures that if the mean-field process starts from the idle initial state, i.e., if $\mf x(0)=\mf e_0=(1,0,0,\dots)$, then it is monotonically non-decreasing in time, i.e.,    
%
\begin{equation}
\label{eqn:g_bound_mf_monotone}
\mf x(t_1, \mf e_0) \leq \mf x(t_2, \mf e_0), \ 0\leq t_1\leq t_2<\infty.
\end{equation}
This follows because the state $\mf e_0$ is dominated by any other state in $S$. In particular, $\mf e_0\leq \mf x(t_2-t_1,\mf e_0)$. Hence, by quasi-monotonicity, we have $\mf x(t_1,\mf e_0)\leq \mf x(t_1, \mf x(t_2-t_1,\mf e_0))=\mf x(t_2,\mf e_0)$.

Furthermore, Lemma~\ref{lem:g_bounded_ODE_Solution_Monotone}
guarantees that if $\mf x(0) \in S$ then $\mf x(t) \in S$ for all $t \geq 0$. Hence, by adding~\eqref{eqn:g_bounded_fluid_process} for all $i \geq k$ and using the fact that $\norm{x(t)}_1 <\infty$ gives
\begin{multline}
\label{eqn:v_n_g_bounded}
 \dot v_k(t)= \lambda \Big[2\epsilon \Big\{x_{k-1}(t) x_{k-1-g}(t)-
 \sum_{i=k}^{k+g-1}x_i(t)\Big(x_{i-1-g}(t)-x_{i-g}(t)\Big) \Big\}+ (1-2\epsilon)(x_{k-1}(t))^2  \Big] -x_{k}(t).
\end{multline}
Since $x_{i-1-g}(t)\geq x_{i-g}(t)$
we have from the above
\begin{align}
\dot v_k&\leq \lambda \Big[2\epsilon x_{k-1} x_{k-1-g}
 + (1-2\epsilon)(x_{k-1})^2  \Big] -x_{k}\nonumber\\
 &= \lambda \Big[\epsilon x_{k-1-g}^2+(1-\epsilon)x_{k-1}^2-\epsilon(x_{k-1-g}-x_{k-1})^2  \Big] -x_{k}\nonumber\\
 &\leq \lambda \Big[\epsilon x_{k-1-g}^2+(1-\epsilon)x_{k-1}^2  \Big] -x_{k}\nonumber\\
 &\leq \lambda x_{k-1-g}^2 -x_k,\label{eq:bound}
\end{align}
where the last inequality follows by using $x_{k-1-g}\geq x_{k-1}$.

{\bf Existence of the Fixed Point $\mf x^*$}:
To prove the existence of the fixed point $\mf x^*$  satisfying~\eqref{eqn:fx_gb_1} and~\eqref{eqn:fx_gb_2}, we first show that $\norm{x_i(t,\mf e_0)}_1$ remains uniformly bounded for all $t\geq 0$. Note that this is a stronger result than $\norm{x_i(t,\mf e_0)}_1 < \infty$ for each $t \geq 0$ which has already been established in Lemma~\ref{lem:g_bounded_ODE_Solution_Monotone}. 

\begin{proposition}
\label{prop:uni_bound}
For $g\geq 0$, $\lambda<1$, let $\mf z^*\in S$ be defined as $z_i^*=1$ for all $i\leq0$ and $z_i^*=\lambda (z_{i-1-g}^*)^2$ for all $i\geq 1$. Then, we have $\mf x(t,\mf e_0) \leq \mf z^*$ for all $t\geq 0$ which implies that
\begin{equation}
\label{eqn:MF_g_bound}
\norm{\mf x(t, \mf e_0)}_1\leq \norm{\mf z^*}_1=(g+1)\sum_{i\geq 1}\lambda^{2^i-1}, \ \forall t\geq0.
\end{equation}
\end{proposition}
\begin{proof}
Let $\mf x(0)=\mf e_0 \leq \mf z^*$. Then, from~\eqref{eqn:g_bound_mf_monotone} we have $\Dot{x}_k(t) \geq 0$ for all $k\geq1$ and for all $t\geq0$. This further implies that $\Dot{v}_k(t)=\sum_{i \geq k}\dot x_i(t) \geq0,\forall k\geq1,\forall t\geq0$. Since $\mf x(0) \leq \mf z^*$, for $\mf x(t) > \mf z^*$ for some $t$ there must exist $t_* < t$ such that $x_l(t_*)=z_l^*$ and $\dot x_l(t_*) > 0$ for some $l\geq1$. Let $m$ be the smallest component where the above two conditions are satisfied. Then, substituting $k=m$ in~\eqref{eq:bound}, we obtain 
\begin{align*}
\Dot{v}_m(t_*) \leq \lambda x_{m-1-g}^2(t_*) -x_m(t_*)&= \lambda x_{m-1-g}^2(t_*) -z_m^*\\
&=\lambda (x_{m-1-g}^2(t_*) - (z_{m-1-g}^*)^2)\\
&\leq 0,
\end{align*}
where the last inequality follows from the definition of $m$. Since we already know that $\dot v_m(t_*) \geq 0$, the above inequality implies $\Dot{v}_m(t_*)=0$. Hence, $\Dot{x}_m(t_*)=\Dot{v}_m(t_*)-\Dot{v}_{m+1}(t_*) = -\Dot{v}_{m+1}(t_*)\leq0$ which contradicts the fact that $\dot x_m(t_*) >0$. 
\end{proof}

 \begin{lem}
 \label{lem:g_bounded_starting_from_0}
Given $g\geq0$, $\lambda<1$, there exists $\mf x^* \in S$ such that
\begin{equation}
\norm{\mf x(t,\mf e_0) -\mf x^{*}}_1 \to 0, \ t \to \infty,
\end{equation}
and $\mf G(\mf x^{*})=\mf 0$. Furthermore, $\mf x^*$ satisfies \eqref{eqn:fx_gb_1} and~\eqref{eqn:fx_gb_2}.
\end{lem}
 \begin{proof}
Since $x_i(t,\mf e_0) \in[0,1]$ for each $i$ and all $t\geq 0$ and $x_i(t,\mf e_0)$  is monotonically non-decreasing in time we must have $x_i(t) \to x_i^{*}$ as $t \to \infty$ for each $i\geq 1$ for some $\mf x^*=(x^*_i) \in \bar S$. 

We first show that the component-wise limit $x^*$ defined above is also the $\ell_1$ limit of $\mf x(t,\mf e_0)$ which will also imply that $\mf x^*\in S$. To show this, we note from Proposition~\ref{prop:uni_bound} that the uniform bound on $\sum_{i \geq 1} x_i(t,\mf e_0)$ in~\eqref{eqn:MF_g_bound} implies by dominated convergence theorem
that 
\begin{equation*}
\lim_{t\to \infty}\norm{\mf x(t)}_1= \sum_{i\geq 1} \lim_{t\to \infty}x_i(t) =\sum_{i\geq 1} x_i^{*}=\norm{\mf x^{*}}_1 \leq (g+1)\sum_{i \geq} \lambda ^{2^i-1}.
\end{equation*}
This shows that $\norm{\mf x(t, \mf e_0)-\mf x^{*}}_1 \to 0$ as $t\to \infty$, and $\mf x^{*} \in S$.

It now remains to show that $\mf G(\mf x^*)=0$. Note that the convergence of $\mf x(t) \to \mf x^{*}$ in $\ell_1$ as $t\to \infty$, and the monotonicity of $\mf x(t)$ imply that for any $\delta>0$ there exists a $t_{\delta}>0$ such that for all $t\geq t_{\delta}$ we have
\begin{multline*}
\delta \geq \norm{\mf x(t+h) - \mf x(t)}_1\geq x_i(t+h) - x_i(t)= \int_t^{t+h} G_i(\mf x(s)) ds\\ \geq h  G_i(\mf x(t^*_{h})), \forall i\geq 1, \forall h\geq0,
\end{multline*}
where  $t^*_{h}\in [t,t+h]$ is the time at which the continuous function $G_i (\mf x(s))$ attains its minimum value in the compact interval $[t,t+h]$. Therefore, we have 
\begin{equation}
\label{eqn:g_bounded_temp1}
 G_i(\mf x(t^*_{h})) \leq \frac{\delta}{h}, \ i\geq1.
\end{equation}
Now we can write
\begin{align}
G_i(\mf x^{*})&=G_i(\mf x^{*}) -G_i(\mf x(t^*_h)) +G_i(\mf x(t^*_h)) \nonumber\\
&\leq \norm{\mf G(\mf x^{*}) - \mf G(\mf x(t^*_h))}_1 +G_i(\mf x(t^*_h))\\
&\leq L_{\lambda}^{\epsilon} \delta  +\frac{\delta}{h},
\end{align}
where for the second inequality we use and the fact that the function $\mf G$ is Lipschitz with constant $L_{\lambda}^{\epsilon}$ and~\eqref{eqn:g_bounded_temp1}. Note that the above inequality is true for any $\delta>0$. Therefore, by fixing $h>0$ and letting $\delta \to 0$ we have $G_i(\mf x^{*})=0$ for all $i\geq1$. Hence, $\mf G (\mf x^{*})=0$.    
Finally, we obtain~\eqref{eqn:fx_gb_1} 
by using $\sum_{i\geq 1} G_i(\mf x^{*})=0$ and~\eqref{eqn:fx_gb_2}
by using $\sum_{i\geq k} G_i(\mf x^{*})=0$.
 \end{proof}

{\bf Global Stability and Uniqueness of the Fixed Point $\mf x^{*}$}:
Now we prove that for any $\mf u \in S$, $\mf x(t,\mf u)$ converges to $\mf x^{*}$ as $t \to \infty$ in $\ell_1$, where $x^*$ is the limit of $\mf x(t, \mf e_0)$ as defined in Lemma~\ref{lem:g_bounded_starting_from_0}.
By Proposition~\ref{prop:uni_bound} and the dominated convergence theorem, it suffices to establish this convergence component-wise. 
Furthermore, it is sufficient to consider
initial points $\mf u \leq \mf x^*$ and $\mf u \geq x^*$ since, by the quasi-monotonicity of $\mf x$, we have
$\mf x(t, \min(\mf u,\mf x^*)) \leq \mf x(t, \mf u) \leq \mf x(t, \max(\mf u,\mf x^*)) $, where the $\min$ and the $\max$ are taken component-wise.

Consider the case when $\mf x(0)=\mf u \leq \mf x^{*}$. Since $\mathbf{u} \geq \mf e_0$, by the quasi-monotonicity of $\mf x$, we have $\mf x(t,\mf e_0)\leq \mf x(t,\mf u) \leq \mf x^{*}, \forall t\geq0.$
%
%
Hence, by Lemma~\ref{lem:g_bounded_starting_from_0}, we have $\mf x(t,\mf u) \to \mf x^*$ since $\mf x(t,\mf e_0) \to \mf x^*$.

Next we consider the case where $\mf x(0)=\mf u \geq \mf x^{*}$. 
We first show that $v_k(t,\mf u)$ remains uniformly bounded for all $t\geq 0$ and for all $k \geq 1$. 
From quasi-monotonicity of $\mf x$ it follows that
$\mf x(t,\mf u) \geq \mf x^{*}$ for all $t \geq 0$. This implies, in particular, that $x_{1}(t,\mf u)\geq x_{1}^{*}=\lambda$.
Hence, from~\eqref{eqn:v_n_g_bounded} for $k=1$, we have $\dot v_1(t,\mf x(0))=\lambda -x_1(t,\mf u) \leq 0$. This implies that
$0\leq v_k(t,\mf u)\leq v_1(t,\mf u)\leq v_1(\mf u)$ for all $t\geq0$ and all $k \geq 1$. 

We shall now establish the convergence $x_{i}(t,\mf u) \to x_{i}^*$  for all $i\geq 1$ by showing 
\eqn{
\label{eqn:int_finite_1}
\int_0^{\infty} (x_{i}(t,\mf u) - x_{i}^*) dt < C_i, 
}
where $C_i >0$ is a finite constant for each $i \geq 1$. To prove~\eqref{eqn:int_finite_1}, we use induction on $i$. For $i=1$, using~\eqref{eqn:fx_gb_1} we have
\eq{
\begin{aligned}
\int_0^{\tau} 
(x_{1}(t,\mf u) - x_{1}^*) dt&=\int_0^{\tau} 
(x_{1}(t,\mf u) -\lambda)=v_1(\mf u)-v_1(\tau,\mf u)\leq v_1(\mf u),
\end{aligned}
}
where the second equality follows from~\eqref{eqn:v_n_g_bounded} for $k=1$ and the inequality follows as $v_1(t,\mf u)$ is uniformly bounded in $t$.
Since the RHS is independent of $\tau$, the integral on the left hand side must be by bounded $v_1(\mf u)$ as $\tau \to \infty$. This shows the base case of the induction. 

Now assume that~\eqref{eqn:int_finite_1} is true for all $i\leq L-1$. For $i=L$, using~\eqref{eqn:v_n_g_bounded} and~\eqref{eqn:fx_gb_2}, we have
\begin{multline}
\int_0^{\tau} 
(x_{L}(t) - x_{L}^{*}) dt\leq v_{L}(\mf u) +  \lambda (1-2\epsilon) \int_0^{\tau}(x_{L-1}^2(t) - (x_{L-1}^{*})^2)dt \\+ \lambda 2\epsilon \int_0^{\tau}(x_{L-1}(t) x_{L-1-g}(t) - x_{L-1}^{*} x_{L-1-g}^{*} )dt- 2\epsilon\lambda \\
\times \int_0^{\tau}\sum_{i=L}^{L+g-1}\sbrac{ x_i(t)(x_{i-1-g}(t)-x_{i-g}(t)) -x_i^{*}(x_{i-1-g}^{*}-x_{i-g}^{*})    }
dt,\nonumber
\end{multline}
where the inequality follows from the uniform boundedness of $v_1(t,\mf u)$ in $t$. To complete the proof, we shall now bound each integral term appearing on the RHS. By using the induction hypothesis and the inequalities $a^2-b^2 \leq 2(a-b)$ and $ab-cd \leq (a-c)+(b-d)$ for $1\geq a\geq c\geq 0$ and $1\geq b \geq d \geq 0$, the integrals in the second and the third terms on the RHS can be easily bounded by $2C_{L-1}$ and $C_{L-1}+C_{L-1-g}$, respectively.

It now remains to bound the integral in the last term. Note that
the third term contains $x_i(t)$
for $i\in \{L,L+1,\ldots, L+g-1\}$
for which the induction hypothesis does not apply. Hence, to bound the integral we need to bound these terms.
We note that by monotonicity we have
$x_i(t) \geq x_i^*$ for all $i\geq 1$.
Hence, $x_i(t)(x_{i-1-g}(t)-x_{i-g}(t)) -x_i^{*}(x_{i-1-g}^{*}-x_{i-g}^{*})\geq 
x_i^*[(x_{i-1-g}(t)-x_{i-1-g}^*)-(x_{i-g}(t)-x_{i-g}^*)]$. Hence, the last term
can be bounded above by
\begin{equation}
\label{eqn:g_bounded_GS_2}
- \lambda 2\epsilon \sum_{i=L}^{L+g-1}x_i^{*}\left[ \int_0^{\tau}(x_{i-1-g}(t)-x_{i-1-g}^{*})dt -\int_0^{\tau}(x_{i-g}(t)-x_{i-g}^{*})dt    \right]\nonumber
.
\end{equation}
Using the induction hypothesis,
we can further bound this by
$2\epsilon \lambda \sum_{i=L}^{L+g-1}(C_{i-1-g}+C_{i-g})$.
This completes the proof of global stability of $\mf x^*$. Since all the trajectories converge to $\mf x^*$, it must be the unique solution of $\mf G(\mf s)=\mf 0$ since starting from any other $\mf y\neq \mf x^*$, satisfying  $\mf G(\mf y)=0$, the trajectory remains at $\mf y$ which contradicts the global stability of $\mf x^*$.

{\bf Limit Interchange}:
Note that Theorem~\ref{thm:stability}.(\ref{stability_g_bounded}), implies that $\pi_n(S)=1,\forall n$. Therefore, we have $\pi_n(\bar{S})=1$ for all $n$.
Since, the space $\bar{S}$ is compact, by Prohorov's theorem the sequence $(\pi_n)_{n}$ must converge weakly to the limit $\pi^*$ with $\pi^*(\bar{S})=1$.
Furthermore, since by~\eqref{eqn:uniform_bound_g_bounded_1}, $\E_{\pi_n} \sbrac{\sum_{i\geq1}x^{(n)}_i(\infty)}$ is uniformly bounded in $n$, we have $\pi^*(S)=1$. 
%
%
Now we prove that the measure $\pi^*$ is the stationary measure of the mean-field process $\mf x$ defined in~\eqref{eqn:g_bounded_fluid_process}.
We know that $(\pi_n)_n \Rightarrow \pi^*$ and the space $S$ is separable. Therefore, the Skorokhod’s Representation 
Theorem implies that $\mf x^{(n)}(0) \overset{a.s}{\to} \mf x(0)$. 
Moreover, if we start the process $\mf x^{(n)}(0)\sim \pi_n$, then $\mf x^{(n)}(t)\sim \pi_n$ for all $t\geq 0$. Hence, from Theorem~\ref{thm:main_g_bounded}.(\ref{po2_g_bounded_Process_level_convergence}) it follows that $\mf x(t) \sim \pi^*$ for all $t\geq0$.
This proves that $\pi^*$ is indeed the stationary measure for the mean-field process $\mf x$. Now from the global stability of the fixed point $\mf x^*,$ it follows immediately that the stationary measure $\pi^*$ is unique and is equal to $\delta_{\mf x^*}$. This completes the proof of limit interchange.

\subsection{Heavy Traffic Limit for the Po2-$(g,\epsilon)$ scheme}
In last we prove Theorem~\ref{thm:main_g_bounded}.(\ref{heavy_gb}). Note that from Lemma~\ref{lem:g_bounded_starting_from_0}, we know $\mf x(t,\mf e_0)$ converges to $\mf x^*$ as $t\to \infty$ in $\ell_1$ and  from Proposition~\ref{prop:uni_bound} we have $$\norm{\mf x(t,\mf e_0)}_1\leq(g+1)\sum_{i\geq 1}\lambda^{2^i-1},$$ for all $t\geq0$. Therefore, we can write $\norm{\mf x^*}_1 \leq (g+1)\sum_{i\geq 1}\lambda^{2^i-1}$.
Hence, dividing the above inequality both side with $-\log(1-\lambda)$, taking limit $\lambda \to 1^-$, we obtain the bound given in~\eqref{eqn:perf_gap_g_bounded}.

\section{Mean-Field Analysis of the Po2-$\epsilon$ Scheme}
\label{sec:po2_eps_analysis}
In this section we prove the main result for the Po2-$\epsilon$ scheme stated in Theorem~\ref{thm:Main_thm_epsilon}.

\subsection{Mean-Field Limit of the  Po2-$\epsilon$ Scheme}
\label{sec:MF_po2_ep}
We first establish the mean-field limit of the Po2-$\epsilon$ policy given in  Theorem~\ref{thm:Main_thm_epsilon}-\ref{po2_eps_Process_level_convergence}. First note that the rate of transitions of the process $\mf x^{(n)}$ from $\mf x \in S^{(n)}$ to $ \mf y \in S^{(n)}$ is given by
\EQN{
\label{eqn:rate_transition}
q^{(n)}_{\mf x, \mf y} =
\begin{cases}
                                      n \lambda p_{i-1}(\mf x),& \text{if $\mf y=\mf x + \mf e_i/n$} \\
                                   n(x_i-x_{i+1}),& \text{if $\mf y= \mf x -\mf e_i/n$}
\end{cases}, \ \forall i\geq1,
}
where $p_{i-1}(\mf x)=(1-\epsilon)(x_{i-1}^2 -x_i^2) + \epsilon ( (1-x_i)^2 - (1-x_{i-1})^2 )$ is the probability that an arrival joins a server of queue length $i-1$.
The mean-field limit of the Po2-$\epsilon$ scheme is proved using the similar argument as shown for the Po2-$(g,\epsilon)$ scheme. First we show that the function $\mf F(\mf x)$ for $\mf x \in S$ defined in~\eqref{eqn:po2_eps_fluid_process} 
is Lipschitz under $\ell_1$-norm. 
\begin{lem}
\label{lem:lip}
The function $\mf F(\mf x)$ is Lipschitz with constant $4\lambda +2$. 
\end{lem}
\begin{proof}
 Let $\mf x, \mf y \in S$. Then we have  
 \begin{align*}
 \norm{\mf F(\mf x) -\mf F(\mf y)}_1&=\sum_{i\geq 1} |  \lambda p_{i-1}(\mathbf{x}) - (x_{i}-x_{i+1}) -  \lambda p_{i-1}(\mathbf{y}) + (y_{i}-y_{i+1}) |\\
 &\leq \lambda\sum_{i\geq 1}|p_{i-1}(\mathbf{x})-p_{i-1}(\mathbf{y})| + 2 \norm{ \mf x -\mf y}_1\\
 &\leq 4\lambda \norm{ \mf x -\mf y}_1 + 2 \norm{ \mf x -\mf y}_1\\
 &=(4\lambda+2)\norm{ \mf x -\mf y}_1.
 \end{align*}
 This completes the proof.
\end{proof}
Furthermore, from~\eqref{eqn:rate_transition} it is clear that the rate at which the jumps occur in $\mf x^{(n)}$ is bounded everywhere that is
\eqn{
\label{eqn:bounded_rate}
\sum_{\mf y \in S}q^{(n)}_{\mf x, \mf y} < n(\lambda+1), \ \forall \ \mf x \in S.
}
Therefore, using Lemma~\ref{lem:lip} and from~\eqref{eqn:bounded_rate} we conclude that the conditions of Kurtz's Theorem are satisfied. Hence, we have 
\begin{equation*}
\lim_{n \to \infty} \sup_{t\leq u} \norm{\mf x^{(n)}(u) - \mf x(u)}_1=0, \ a.s.
\end{equation*}
This completes the proof.

\subsection{Mean-Field Steady State Behaviour for the Po2-$\epsilon$ Scheme}
\label{sec:MF_eps_prop}
In this section we
prove Theorem~\ref{thm:Main_thm_epsilon}.(\ref{fixed_point}), which shows that the differential equations defined in~\eqref{eqn:po2_eps_fluid_process} has a unique fixed point $\mf x^*$ and
it follows the recursion defined in~\eqref{eqn:po2_eps_Fixed_Point}. Moreover, we prove that the fixed point $\mf x^*$ is globally stable and finally establish the interchange of limits. 

For $\mf u \in S$, we define $v_{k}(t,\mf u)=\sum_{i\geq k}x_{i}(t,\mf u)$ and $v_{k
}(\mf u)=\sum_{i\geq k}u_{i}$
for each $k\geq1$. 

\begin{lem}
\label{lem:po2_eps_ODE_Solution_Monotone}
Let $\epsilon\in [0,1]$. The following statements hold for the process $\mf x$ defined in Theorem~\ref{thm:Main_thm_epsilon}. 
\begin{enumerate}
\item If $\mathbf{u}\in S$, then, for any $t\geq 0$, we have $\mf x(t,\mf u) \in S$.
\item For any $\mf u, \mf u' \in S$ satisfying $\mf u \leq \mf u'$
we have ${\mf x}(t, \mf u) \leq {\mf x}(t, \mf u')$ for all $t \geq 0$, where the inequality $\leq$ is understood component-wise.
\end{enumerate}
\end{lem}
\begin{proof}
From~\eqref{eqn:po2_eps_fluid_process} and~\eqref{eqn:arrival_prop_eps}, we can write
\begin{align}
 \dot x_i \leq \lambda p_{i-1}(\mf x) &\leq \lambda [(1-2\epsilon)(x^2_{i-1}-x_i^2)+2\epsilon (x_{i-1}-x_i)]\nonumber\\
 &\leq \lambda \max(1,2\epsilon)x_{i-1}
=M_{\lambda,\epsilon}x_{i-1},\nonumber
\end{align}
where $M_{\lambda,\epsilon}=\lambda\max(1,2\epsilon)$.
Therefore, for each $i \geq 1$ we have $$x_i(t) \leq x_i(0)+M_{\lambda,\epsilon}\int_{0}^t x_{i-1}(s) ds.$$
Hence, using the above recursively for each $i \geq 1$, we obtain
\begin{equation*}
x_i(t) \leq x_i(0)+\sum_{k=0}^{i-1}x_{k}(0)\frac{(M_{\lambda,\epsilon}t)^{i-k}}{(i-k)!}.
\end{equation*}
Summing the above for all $i \geq 1$
we obtain $v_1(\mf x(t))\leq (1+v_1(\mf x(0)))\exp(M_{\lambda,\epsilon}t).$
Therefore, if $v_1(\mf x(0))<\infty$, then $v_1(\mf x(t))<\infty$ for all $t$. This completes the proof of first part.

To prove second part we need to show that $\frac{d{x}_{i}(t)}{dt}$ is non-decreasing in $x_j(t)$ for all $j\neq i$. We know from~\eqref{eqn:po2_eps_fluid_process} that 
\EQ{
\frac{d{x}_{i}(t)}{dt}= \lambda \sbrac{(1-2\epsilon)(x_{i-1}^2(t)- x_i^2(t)) + 2\epsilon (x_{i-1}(t) -x_i(t))}  -(x_{i}(t)-x_{i+1}(t)).
}
Therefore, it is clear that the above expression is non-decreasing with $x_{i+1}(t)$. Note that the derivative of the terms involving $x_{i-1}(t)$ component in the above expression is $2(1-2\epsilon)x_{i-1}(t) +2\epsilon$, which clearly positive for $\epsilon\leq 1/2$. Moreover, for $\epsilon>1/2$ we can write the derivative of terms involving $x_{i-1}(t)$ as $2[\epsilon - (2\epsilon -1) x_{i-1}(t)]$ which is decreasing with $x_{i-1}(t)$ and has minimum value of $2(1-\epsilon)>0$. Hence, $\frac{d{x}_{i}(t)}{dt}$ is also non-decreasing with $x_{i-1}(t)$.
 \end{proof}
{\bf Fixed Point}: 
For fixed point $\mf x^*$ we need to equate~\eqref{eqn:po2_eps_fluid_process} to $0$ and get
\begin{equation}
\label{eqn:fix_temp}
\lambda \sbrac{(1-2\epsilon)(x_{i-1}^{*2}- x_i^{*2}) + 2\epsilon (x_{i-1}^* -x_i^*)} =x_{i}^*-x_{i+1}^*, \ i\geq 1.  \end{equation}
Summing~\eqref{eqn:fix_temp} for all $i\geq 1$ we get $x_1^*=\lambda$. Moreover, summing~\eqref{eqn:fix_temp} for all $i\geq j$ we obtain
\begin{equation*}
x^*_{j}=\lambda \sbrac{(1-2\epsilon)(x_{j-1}^*)^2 + 2\epsilon x_{j-1}^*}, \ \forall j\geq2.
\end{equation*}

{\bf Global Stability}:
To prove global stability of the fixed point $\mf x^*$, we use the monotonicity of the mean-field
process $\mf x$ shown in Lemma~\ref{lem:po2_eps_ODE_Solution_Monotone}.

Note that Lemma~\ref{lem:po2_eps_ODE_Solution_Monotone}
guarantees that if $\mf x(0) \in S$ then $\mf x(t) \in S$ for all $t \geq 0$. Hence, by adding~\eqref{eqn:po2_eps_fluid_process} for all $i \geq k$ and using the fact that $\norm{x(t)}_1 <\infty$ gives
\begin{equation}
\label{eqn:v_n}
 \frac{dv_{k}(t,\mf u)}{dt}=   \lambda \sbrac{(1-2\epsilon)x_{k-1}^2(t,\mf u) + 2\epsilon x_{k-1}(t,\mf u) }-x_{k}(t,\mf u)
\end{equation}
Specifically, for $k=1$ we have 
\begin{equation}
\label{eqn:v_1}
  \frac{dv_{1}(t,\mf u)}{dt}= \lambda -x_1(t,\mf u). 
\end{equation}
Now from the monotonicity property (Lemma~\ref{lem:po2_eps_ODE_Solution_Monotone}) of the mean-field process $\mf x$ we have for any $\mf x(0) \in S$
\EQN{
\label{eqn:sandwitch_sol}
\mathbf{x}(t,\min(\mathbf{x}(0),\mathbf{x}^*))\leq\mathbf{x}(t,\mathbf{x}(0))\leq \mathbf{x}(t,\max(\mathbf{x}(0),\mathbf{x}^*)), \ t\geq 0,
}
where $\min(\mf u,\mf v)$ with $\mf u, \mf v \in S$
is defined by taking the component-wise minimum. From~\eqref{eqn:sandwitch_sol} it is clear that to prove global stability it is enough to prove convergence $\mf x(t,\mf x(0)) \to \mf x^*$ holds for
initial states satisfying either of the following two conditions:
(i) $\mf x(0) \geq \mf x^*$ and (ii) $\mf x(0) \leq \mf x^*$.

To prove convergence holds for above two initial conditions we first show that for any solution $\mf x(\cdot,\mf x(0)) \in S$, $v_k(t,\mf x(0))$ is uniformly bounded in $t$ for all $k \geq 1$.
Consider the case when $\mf x(0) \geq \mf x^*$. 
From Lemma~\ref{lem:po2_eps_ODE_Solution_Monotone}
it follows that for $\mf x(0) \geq \mf x^*$, we have
$\mf x(t,\mf x(0)) \geq \mf x^*$ for all $t \geq 0$. Therefore, we can write
\eq{
 x_{1}(t,\mf x(0))\geq x_{1}^*=\lambda,
}
where the last equality follows from~\eqref{eqn:po2_eps_Fixed_Point}.
Hence, from~\eqref{eqn:v_1} we have $\frac{dv_1(t,\mf x(0))}{dt}\leq 0$
from which it follows that
$0\leq v_1(t,\mf x(0))\leq v_1(\mf x(0))$ for all $t\geq0$. 
Since the sequence $\brac{v_k(t,\mf x(0))}_{k\geq 1}$ is non-increasing, we have $0\leq v_k(t,\mf x(0))\leq v_1(\mf x(0))$ for all $k\geq1$ and for all $t\geq0$. This proves that $v_k(t,\mf x(0))$ is uniformly bounded in $t$ for each $k\geq1$ if $\mf x(0) \geq \mf x^*$.
Now consider the case $\mf x(0)\leq \mf x^*$. From Lemma~\ref{lem:po2_eps_ODE_Solution_Monotone}
it follows that for $\mf x(0) \leq \mf x^*$, we have
$\mf x(t,\mf x(0)) \leq \mf x^*$ for all $t \geq 0$. Therefore, we have $v_1(t,\mf x(0))\leq v_1(\mf x^*)$ for all $t\geq0$. This shows that the component $v_k(t,\mf x(0))$ is uniformly bounded in $t$ for each $k\geq1$ for $\mf x(0)\leq \mf x^*$.

Since  $v_k(t,\mf x(0))$ is uniformly bounded in $t$, the convergence $x_{i}(t,\mf x(0)) \to x_{i}^*$  for all $i\geq 1$ will follow from 
\eqn{
\label{eqn:int_finite}
\int_0^{\infty} (x_{i}(t,\mf x(0)) - x_{i}^*) dt < \infty, \ \forall i\geq1,
}
for the case $\mf x(0)\geq \mf x^*$ and from
\eqn{
\label{eqn:int_finite_2}
\int_0^{\infty} (x_{i}^* - x_{i}(t,\mf x(0))) dt < \infty, \ i\geq1,
}
for the case $\mf x(0)\leq \mf x^*$. We now prove~\eqref{eqn:int_finite} to show convergence for the case $\mf x(0)\geq \mf x^*$; the proof of other case follows similarly.

We will use induction starting with $i=1$. We can write~\eqref{eqn:int_finite} for $i=1$ as
\eq{
\begin{aligned}
\int_0^{\tau} 
(x_{1}(t,\mf x(0)) - x_{1}^*) dt&=\int_0^{\tau} 
(x_{1}(t,\mf x(0)) -\lambda)\\
&=-  \int_0^{\tau} \frac{dv_1(t,\mf x(0))}{dt}dt\\
&=v_1(\mf x(0))-v_1(\tau,\mf x(0))\\
&\leq v_1(\mf x(0)),
\end{aligned}
}
where the second equality follows from~\eqref{eqn:v_1} and the inequality follows as $v_1(t,\mf x(0))$ is uniformly bounded in $t$.
Observe that the right hand side is bounded by a constant for all $\tau$ , the integral on the left hand side must converge as $\tau \to \infty$. This shows that $x_{1}(t,\mf x(0)) \to x_{1}^*$ as $t\to \infty$.
Now assume that~\eqref{eqn:int_finite} is true for all $i\leq L-1$. For $i=L$ we can write~\eqref{eqn:int_finite} as
\begin{equation*}
\begin{aligned}
 \int_0^{\tau} 
(x_{L}(t,\mf x(0)) - x_{L}^*) dt&= \int_0^{\tau} \Big[ \frac{-dv_L(t,\mf x(0))}{dt} + \lambda \cbrac{(1-2 \epsilon)x_{L-1}^2(t,\mf x(0)) + 2\epsilon x_{L-1}(t,\mf x(0)) }\\
&-\lambda (1-2\epsilon)(x_{L-1}^*)^2 -\lambda2\epsilon x_{L-1}^* \Big]dt,\\
&=v_L(\mf x(0))-v_L(\tau,\mf x(0)) + \int_0^{\tau} \Big[ \lambda(1-2\epsilon) ( x_{L-1}^2(t,\mf x(0))-(x_{L-1}^*)^2  )\\
&+ \lambda2\epsilon (x_{L-1}(t,\mf x(0)) -x_{L-1}^*) \Big]dt,\\
&\leq v_L(\mf x(0)) + \int_0^{\tau}\lambda \Big[ (1-2\epsilon) ( x_{L-1}^2(t,\mf x(0))-(x_{L-1}^*)^2  )\Big]dt\\
&+ \int_0^{\tau}\lambda \Big[2\epsilon (x_{L-1}(t,\mf x(0)) -x_{L-1}^*) \Big]dt 
\end{aligned}
\end{equation*}
where the first equality follows from~\eqref{eqn:v_n} for $k=L$ and from~\eqref{eqn:po2_eps_Fixed_Point} for $i=L$. Moreover, the inequality follows as $v_L(t,\mf x(0))$ is uniformly bounded in $t$.
Note that by the induction hypothesis, the last two integral on the right hand side of above expression converges as $\tau \to \infty$.
Hence, the
integral on the left hand side also must converge as required.

{\bf Limit Interchange}:
We know from Theorem~\ref{thm:stability}.(\ref{epsilon_stability}) that for $\epsilon\in[0,1]$, and $\lambda<\min(1,1/2\epsilon)$ we have 
$$\E_{\pi_n} \sbrac{\sum_{i\geq1}x^{n}_i(\infty)}\leq  \frac{\lambda}{1-\max(1,2\epsilon)\lambda}.$$
Therefore, using the global stability result and the process convergence result of the Po2-$\epsilon$ scheme, the limit interchange follows using similar argument as proved for the Po2-$(g,\epsilon)$ scheme.

\subsection{Heavy-Traffic Limit of the Po2-$\epsilon$ Scheme}
In this section we prove  Theorem~\ref{thm:Main_thm_epsilon}.(\ref{comp_ratio_epsilon}), which computes the ratio of the average response time of jobs under the Po2-$\epsilon$ scheme with the logarithmic of average response time of jobs under the random scheme as $\lambda \to 1$. We first write the recursion given in~\eqref{eqn:po2_eps_Fixed_Point} as 
\begin{equation}
\label{eqn:recursion}
 x_{i}^* = T_{\lambda,\epsilon}(x_{i-1} ^*)= \lambda \sum_{i=1}^2 a_i(x_{i-1}^*)^{d_i}, \ \forall i\geq 1,
\end{equation}
where $a_1=(1-2\epsilon)$, $a_2=2\epsilon$, $d_1=2$, and $d_2=1$. 
The result now follows from Theorem~4.5 of~\cite{hellemans2020heavy}.

\section{Conclusion and Future Directions}
\label{sec:Conclusion}
In this paper, we analyzed the effects of load comparison errors on the performance of the Po2 scheme. We considered two models of error. For the load-dependent error model, we showed that the Po2 scheme retains its benefits over the random scheme in the heavy traffic limit $\lambda \to 1$ for all values of $g$ and $\epsilon$. 
For the load-independent error model, we have shown that the Po2 scheme retains its benefits over the random scheme only if the probability of error $\epsilon\leq 1/2$. We introduce a general framework using Lyapunov functions to prove stability of our schemes. We also use a new approach to establish the mean-field limit results as the fixed point does not admit a recursive solution. 

There are many interesting directions for further research. We have analyzed the performance of the Po2-$(g,\epsilon)$ scheme assuming $g$ to be constant independent of $n$.  It will be interesting to see the effect of varying $g$ as a function of $n$. Another direction is to study the effects of delay in receiving the queue length information at the dispatcher. A more explicit delay dependent error model can be considered. Here, the challenge will be to analyze the effect of the delay on the performance of the Po2 scheme.

\bibliographystyle{IEEEtran}
\bibliography{sample}

\begin{thebibliography}{10}
\providecommand{\url}[1]{#1}
\csname url@samestyle\endcsname
\providecommand{\newblock}{\relax}
\providecommand{\bibinfo}[2]{#2}
\providecommand{\BIBentrySTDinterwordspacing}{\spaceskip=0pt\relax}
\providecommand{\BIBentryALTinterwordstretchfactor}{4}
\providecommand{\BIBentryALTinterwordspacing}{\spaceskip=\fontdimen2\font plus
\BIBentryALTinterwordstretchfactor\fontdimen3\font minus
  \fontdimen4\font\relax}
\providecommand{\BIBforeignlanguage}[2]{{%
\expandafter\ifx\csname l@#1\endcsname\relax
\typeout{** WARNING: IEEEtran.bst: No hyphenation pattern has been}%
\typeout{** loaded for the language `#1'. Using the pattern for}%
\typeout{** the default language instead.}%
\else
\language=\csname l@#1\endcsname
\fi
#2}}
\providecommand{\BIBdecl}{\relax}
\BIBdecl

\bibitem{Weber1978}
R.~R. Weber, ``On the optimal assignment of customers to parallel servers,''
  \emph{Journal of Applied Probability}, vol.~15, no.~2, pp. 406--413, 1978.

\bibitem{Winston1977optimality}
W.~Winston, ``Optimality of the shortest line discipline,'' \emph{Journal of
  applied probability}, vol.~14, no.~1, pp. 181--189, 1977.

\bibitem{sanidhay_performance}
S.~Bhambay and A.~Mukhopadhyay, ``Asymptotic optimality of speed-aware jsq for
  heterogeneous service systems,'' \emph{Performance Evaluation}, vol. 157, p.
  102320, 2022.

\bibitem{sanidhay_wiopt}
------, ``Optimal load balancing in heterogeneous server systems,'' in
  \emph{2022 20th International Symposium on Modeling and Optimization in
  Mobile, Ad hoc, and Wireless Networks (WiOpt)}.\hskip 1em plus 0.5em minus
  0.4em\relax IEEE, 2022, pp. 113--120.

\bibitem{goren2022distributed}
G.~Goren, S.~Vargaftik, and Y.~Moses, ``Distributed dispatching in the parallel
  server model,'' \emph{IEEE/ACM Transactions on Networking}, 2022.

\bibitem{goren2021stochastic}
------, ``Stochastic coordination in heterogeneous load balancing systems,'' in
  \emph{Proceedings of the 2021 ACM Symposium on Principles of Distributed
  Computing}, 2021, pp. 403--414.

\bibitem{Lu_JIQ_2011}
Y.~Lu, Q.~Xie, G.~Kliot, A.~Geller, J.~R. Larus, and A.~Greenberg,
  ``Join-idle-queue: A novel load balancing algorithm for dynamically scalable
  web services,'' \emph{Perform. Eval.}, vol.~68, no.~11, p. 1056–1071, nov
  2011.

\bibitem{stolyar_JIQ}
A.~L. Stolyar, ``Pull-based load distribution among heterogeneous parallel
  servers: the case of multiple routers,'' \emph{Queueing Systems}, vol.~85,
  no. 1-2, pp. 31--65, 2017.

\bibitem{Mitzenmacherthesis}
M.~Mitzenmacher, ``The power of two choices in randomized load balancing,''
  \emph{PhD thesis, University of California at Berkeley}, 1996.

\bibitem{Vvedenskaya1996}
N.~D. Vvedenskaya, R.~L. Dobrushin, and F.~I. Karpelevich, ``Queueing system
  with selection of the shortest of two queues: An asymptotic approach,''
  \emph{Problemy Peredachi Informatsii}, vol.~32, no.~1, pp. 20--34, 1996.

\bibitem{tarreau2019test}
\BIBentryALTinterwordspacing
W.~Tarreau, ``Test driving power of two random choices load balancing,'' 2019.
  [Online]. Available:
  \url{https://www.haproxy.com/blog/power-of-two-load-balancing/}
\BIBentrySTDinterwordspacing

\bibitem{garrett2018nginx}
\BIBentryALTinterwordspacing
O.~Garrett, ``Nginx and the ‘power of two choices’ load-balancing
  algorithm,'' 2018. [Online]. Available: \url{https://www.nginx.com/blog/}
\BIBentrySTDinterwordspacing

\bibitem{mike2018netflix}
\BIBentryALTinterwordspacing
M.~Smith, ``Rethinking netflix’s edge load balancing,'' 2018. [Online].
  Available:
  \url{https://netflixtechblog.com/netflix-edge-load-balancing-695308b5548c}
\BIBentrySTDinterwordspacing

\bibitem{zhou2021asymptotically}
X.~Zhou, N.~Shroff, and A.~Wierman, ``Asymptotically optimal load balancing in
  large-scale heterogeneous systems with multiple dispatchers,''
  \emph{Performance Evaluation}, vol. 145, pp. Art--No, 2021.

\bibitem{burke2021misreporting}
Q.~Burke, P.~McDaniel, T.~La~Porta, M.~Yu, and T.~He, ``Misreporting attacks
  against load balancers in software-defined networking,'' \emph{Mobile
  networks and applications}, 2021.

\bibitem{sitaraman2001power}
M.~Mitzenmacher, A.~W. Richa, and R.~Sitaraman, ``The power of two random
  choices: A survey of techniques and results,'' \emph{Handbook of randomized
  computing, Netherlands, 255–312}, 2001.

\bibitem{los2022balanced}
D.~Los and T.~Sauerwald, ``Balanced allocations with the choice of noise,'' in
  \emph{Proceedings of the 2022 ACM Symposium on Principles of Distributed
  Computing}, 2022, pp. 164--175.

\bibitem{nadiradze2021achieving}
G.~Nadiradze, ``On achieving scalability through relaxation,'' \emph{PhD
  thesis}, 2021.

\bibitem{bramson_stability}
M.~Bramson, ``Stability of join the shortest queue networks,'' \emph{The Annals
  of Applied Probability}, vol.~21, no.~4, pp. 1568--1625, 2011.

\bibitem{foss_stability}
S.~Foss and N.~Chernova, ``On the stability of a partially accessible
  multi-station queue with state-dependent routing,'' \emph{Queueing Systems},
  vol.~29, pp. 55--73, 1998.

\bibitem{bramson2012asymptotic}
M.~Bramson, Y.~Lu, and B.~Prabhakar, ``Asymptotic independence of queues under
  randomized load balancing,'' \emph{Queueing Systems}, vol.~71, no.~3, pp.
  247--292, 2012.

\bibitem{bramson2010randomized}
------, ``Randomized load balancing with general service time distributions,''
  \emph{ACM SIGMETRICS performance evaluation review}, vol.~38, no.~1, pp.
  275--286, 2010.

\bibitem{chen2012asymptotic}
H.~Chen and H.-Q. Ye, ``Asymptotic optimality of balanced routing,''
  \emph{Operations research}, vol.~60, no.~1, pp. 163--179, 2012.

\bibitem{maguluri2014heavy}
S.~T. Maguluri, R.~Srikant, and L.~Ying, ``Heavy traffic optimal resource
  allocation algorithms for cloud computing clusters,'' \emph{Performance
  Evaluation}, vol.~81, pp. 20--39, 2014.

\bibitem{mukherjee2020asymptotic}
D.~Mukherjee, S.~C. Borst, J.~S. Van~Leeuwaarden, and P.~A. Whiting,
  ``Asymptotic optimality of power-of-d load balancing in large-scale
  systems,'' \emph{Mathematics of Operations Research}, vol.~45, no.~4, pp.
  1535--1571, 2020.

\bibitem{budhiraja2019supermarket}
A.~Budhiraja, D.~Mukherjee, and R.~Wu, ``Supermarket model on graphs,''
  \emph{The Annals of Applied Probability}, vol.~29, no.~3, pp. 1740--1777,
  2019.

\bibitem{debankur_constrained_2021}
D.~Rutten and D.~Mukherjee, ``Load balancing under strict compatibility
  constraints,'' ser. SIGMETRICS '21.\hskip 1em plus 0.5em minus 0.4em\relax
  New York, NY, USA: Association for Computing Machinery, 2021, p. 51–52.

\bibitem{zhao2022exploiting}
Z.~Zhao, D.~Mukherjee, and R.~Wu, ``Exploiting data locality to improve
  performance of heterogeneous server clusters,'' \emph{arXiv preprint
  arXiv:2211.16416}, 2022.

\bibitem{arpan_tcns}
A.~Mukhopadhyay and R.~R. Mazumdar, ``Analysis of randomized
  join-the-shortest-queue (jsq) schemes in large heterogeneous
  processor-sharing systems,'' \emph{IEEE Transactions on Control of Network
  Systems}, vol.~3, no.~2, pp. 116--126, 2015.

\bibitem{arpan_ssy}
A.~Mukhopadhyay, A.~Karthik, and R.~R. Mazumdar, ``Randomized assignment of
  jobs to servers in heterogeneous clusters of shared servers for low delay,''
  \emph{Stochastic Systems}, vol.~6, no.~1, pp. 90--131, 2016.

\bibitem{Kelly_book}
F.~Kelly and E.~Yudovina, \emph{Stochastic Networks}.\hskip 1em plus 0.5em
  minus 0.4em\relax Cambridge University Press, 2014.

\bibitem{Glynn_bounds}
P.~W. Glynn and A.~Zeevi, ``Bounding stationary expectations of markov
  processes,'' in \emph{Markov processes and related topics: a Festschrift for
  Thomas G. Kurtz}.\hskip 1em plus 0.5em minus 0.4em\relax Institute of
  Mathematical Statistics, 2008, pp. 195--214.

\bibitem{srikant2013communication}
R.~Srikant and L.~Ying, \emph{Communication networks: an optimization, control,
  and stochastic networks perspective}.\hskip 1em plus 0.5em minus 0.4em\relax
  Cambridge University Press, 2013.

\bibitem{kurtz_book}
S.~N. Ethier and T.~G. Kurtz, \emph{Markov processes: characterization and
  convergence}.\hskip 1em plus 0.5em minus 0.4em\relax John Wiley \& Sons,
  2009, vol. 282.

\bibitem{kurtz1981approximation}
T.~G. Kurtz, \emph{Approximation of population processes}.\hskip 1em plus 0.5em
  minus 0.4em\relax SIAM, 1981.

\bibitem{Deimling2006}
K.~Deimling, \emph{Ordinary differential equations in Banach spaces}.\hskip 1em
  plus 0.5em minus 0.4em\relax Springer, 2006, vol. 596.

\bibitem{hellemans2020heavy}
T.~Hellemans and B.~Van~Houdt, ``Heavy traffic analysis of the mean response
  time for load balancing policies in the mean field regime,'' \emph{arXiv
  preprint arXiv:2004.00876}, 2020.

\end{thebibliography}

\end{document}